\documentclass[12pt]{article}
\usepackage[latin1]{inputenc}
\usepackage{amsfonts,amssymb,amsmath}
\usepackage{amsthm}
\usepackage{amstext,amsfonts, amscd}
\usepackage{xcolor}
\usepackage{hyperref}          
\usepackage{ulem}

\def\Box{\leavevmode\vbox{\hrule
     \hbox{\vrule\kern4pt\vbox{\kern4pt}%
           \vrule}\hrule}}
\def\blackbox{\leavevmode\vrule height 5pt width 4pt depth 0pt\relax}
\def\endproof{\null\hfill {$\blackbox$}\bigskip}

\def\paragraph#1{{\bf #1\ }}

\newcommand{\rot}{\Theta}

\newtheorem{lemma}{Lemma}[section]  

\newtheorem{theorem}[lemma]{Theorem}

\newtheorem{proposition}[lemma]{Proposition}

\newtheorem{hypothesis}{Hypothesis}[section]

\title{Hyperbolicity and non-conservativity of a hydrodynamic model of swarming rigid bodies} 
\author{P. Degond$^{(1)}$, A. Frouvelle$^{(2)}$, S. Merino-Aceituno$^{(3)}$, A. Trescases$^{(1)}$} 
\date{} 
\begin{document}

\maketitle

\vspace{0.5 cm}

\begin{center}

$^{(1)}$ Institut de Math\'ematiques de Toulouse ; UMR5219 \\
Universit\'e de Toulouse ; CNRS \\
UPS, F-31062 Toulouse Cedex 9, France\\
email PD: pierre.degond@math.univ-toulouse.fr\\
email AT: ariane.trescases@math.univ-toulouse.fr

\bigskip

$^{(2)}$ CEREMADE, CNRS, Universit\'e Paris-Dauphine\\
Universit\'e PSL, 75016 Paris, France \\
email: frouvelle@ceremade.dauphine.fr

\bigskip

$^{(3)}$ Faculty of Mathematics, University of Vienna,\\
Oskar-Morgenstern-Platz 1, 1090 Vienna, Austria\\
email: sara.merino@univie.ac.at

\end{center}

\vspace{0.5 cm}
\begin{abstract}
In this paper, we study a nonlinear system of first order partial differential equations describing the macroscopic behavior of an ensemble of interacting self-propelled rigid bodies. Such system may be relevant for the modelling of bird flocks, fish schools or fleets of drones. We show that the system is hyperbolic and can be approximated by a conservative system through relaxation. We also derive viscous corrections to the model from the hydrodynamic limit of a kinetic model. This analysis prepares the future development of numerical approximations of this system. 
\end{abstract}

\medskip
\noindent
{\bf Acknowledgements:} PD holds a visiting professor association with the Department of Mathematics, Imperial College London, UK.
SM's research was supported by the Austrian Science Fund (FWF) through the  project  F65 and by the Vienna Science  and  Technology  Fund  (WWTF)  with  a  Vienna  Research  Groups  for  Young  Investigators, grant VRG17-014. AF acknowledges support from the Project EFI ANR-17-CE40-0030 of the French National Research Agency. AF thanks the hospitality of the Laboratoire de Math\'ematiques et Applications (LMA, CNRS) in the Universit\'e de Poitiers, where part of this research was conducted.


\vskip 0.4cm

\setcounter{equation}{0}
\section{Introduction}
\label{intro}

In this paper, we study a new hydrodynamic model of swarming behavior derived in~\cite{degondfrouvellemerino17, degondfrouvellemerinotrescases18, degond2018alignment} in dimension 3 and later in~\cite{degond2021body} in arbitrary dimension. This model describes the motion of an ensemble of self-propelled rigid bodies that tend to adjust their body-frame to the averaged body-frame of their neighbors, up to some noise. It may be used to describe e.g. bird flocks, fish schools or fleets of drones. The model consists of a system of coupled first-order partial differential equations in non-conservative form for the agents' local density~$\rho$ and local average frame (represented by a rotation matrix~$\rot$ mapping a fixed reference frame to this frame). We show that this model is hyperbolic and provide two approaches to deal with the issues brought by the non-conservative character of the model. 

Collective dynamics, of which swarming behavior is one of the manifestations, occurs ubiquitously in the living world, from embryonic cell  migration~\cite{giniunaite2020modelling} to human crowds~\cite{warren2018collective} through locust swarms~\cite{bazazi2008collective} or fish schools~\cite{lopez2012behavioural}. Swarming behavior occurs when all agents reach a consensus about their direction of motion. Other manifestations of collective dynamics are, among others, travelling-waves~\cite{welch2001cell}, oscillations~\cite{petrolli2021oscillations} or segregation~\cite{kabla2012collective}. These effects occur at the large (or macroscopic) scale, i.e. scales that are much larger than the typical interaction range between the agents. This motivates the use of macroscopic models to analyze them. 

Macroscopic models describe a large population of agents by means of continuum fields (such as the agents' density or mean velocity) which solve systems of partial differential equations. Examples of  macroscopic models used in collective dynamics can be found in~\cite{bertin2009hydrodynamic, bertozzi2015ring, degond2008continuum, ha2008from, toner1998flocks}. By contrast, microscopic models describe the behavior of each agent individually and lead to large systems of ordinary or stochastic differential equations. They provide a more accurate description of the particle swarm but at the expense of an increased computational complexity. They are also less amenable to theoretical analysis such as stability, asymptotic behavior, etc. Examples of such microscopic models can be found in~\cite{aoki1982simulation, chate2008collective, cucker2007emergent, czirok1996formation, dorsogna2006self, ha2009simple, motsch2011new, vicsek1995novel} (see a review in~\cite{vicsek2012collective}). There is a systematic way to pass from microscopic models to macroscopic ones. This involves an intermediate class of models called kinetic or mean-field models, where the particle dynamics is described statistically in terms of a probability distribution function. This methodology is classically used in rarefied gas dynamics (see e.g.~\cite{cercignani2013mathematical, degond2004macroscopic}). Its usage in collective dynamics can be found in~\cite{degondfrouvelleliu15, degond2008continuum, ha2008from, jiang2016hydrodynamic}. Kinetic models of collective behavior have been studied for their own sake or in relation to particle models in~\cite{bertin2006boltzmann, bolley2012mean, bostan2013asymptotic, briant2020cauchy, frouvelle2012dynamics, ha2009simple, peruani2008mean}. 

The macroscopic model studied in this paper is referred to as the Self-Organized Hydrodynamics for Body-orientation (SOHB). It originates from a microscopic model of interacting self-propelled rigid bodies described in~\cite{degond2021body, degondfrouvellemerino17, degondfrouvellemerinotrescases18, degond2018alignment} and simulated in~\cite{degond2021bulk} which will not be reproduced here. Particle dynamics involving the full rigid-body attitude of the agents have been previously proposed in~\cite{fetecau2021emergent, golse2019mean, ha2017emergent, hildenbrandt2010self}, although in a different form than used here. The corresponding kinetic model will be presented in Section~\ref{sec:viscosity} but most of the paper will be focused on the SOHB model which will be presented in Section~\ref{sec_SOHB_presentation}. The rigorous derivation of the kinetic model has been investigated in~\cite{diez2019propagation}. With a slight modification this kinetic model produces phase transitions which have been studied in~\cite{degond2019phase, frouvelle2021body}. The SOHB model exhibits topologically non-trivial solutions~\cite{degond2021bulk}. 

The present study parallels that of the Vicsek model, where self-propelled particle interact by locally aligning their self-propulsion velocity~\cite{vicsek1995novel}. The corresponding kinetic model has been derived in~\cite{bolley2012mean, briant2020cauchy} and studied theoretically in~\cite{bostan2013asymptotic, figalli2018global, gamba2016global} and numerically in~\cite{dimarcomotsch16, gamba2015spectral, griette2019kinetic}. Phase transitions in the kinetic model have been investigated in~\cite{degond2013macroscopic, degondfrouvelleliu15, frouvelle2012dynamics}. The macroscopic model associated with this particle system has been referred to as the Self-Organized Hydrodynamics (SOH). It has first been derived in~\cite{degond2008continuum} with further elaborations in~\cite{degond20133hydrodynamic, frouvelle2012continuum}. The convergence from the kinetic to the fluid model has been rigorously proved in~\cite{jiang2016hydrodynamic}. Other macroscopic versions of the Vicsek particle system have been proposed in the literature~\cite{toner1998flocks} but no mathematical theory is available yet to support them as approximations of the Vicsek model. Numerical simulations of the SOH model have been provided in~\cite{degond2015macroscopic, dimarcomotsch16, motsch2011numerical}. As this description shows, the mathematical theory of the SOH and related models is more complete than that of the SOHB and related models. The study of the latter still involves many open questions. At several instances in this paper, reference will be made to the SOH model in comparison with the SOHB model. 

The purpose of this paper is to study two aspects of the SOHB model. First, in Section~\ref{sec_hyperbolic} we will show that, as a nonlinear system of first order partial differential equations, it is hyperbolic, a good indication that it is well-posed, at least locally. The proof may be cumbersome, unless an appropriate reference frame is used, and the model shown to be invariant by this choice of frame. By contrast, the hyperbolicity of the SOH model was more straightforward (and shown in~\cite{degond2008continuum}). 

The second aspect relates to the non-conservativity of the SOHB model, i.e. the fact that the spatial derivatives are not in divergence form. It is well-known that nonlinear hyperbolic systems may generate discontinuous weak solutions, referred to as shock waves or contact discontinuities. In the case of conservative systems, these solutions can be unequivocally defined~\cite{godlewski2013numerical, leveque1992numerical, smoller2012shock} thanks to the Rankine-Hugoniot relations (supplemented by entropy conditions in the case of shocks). For non-conservative systems, there is no unique way to define shock solutions. Rather, there exists a whole family of solutions depending on how the so-called non-conservative product can be defined~\cite{lefloch2002hyperbolic}. This multiplicity of solutions is due to loss of information in the passage from the microscopic to macroscopic dynamics. Indeed, dissipation, which occurs at the microscopic scale, eventually decides what expression of the non-conservative product needs to be used. However, these dissipation mechanisms are infinitesimal at the macroscopic scale and are lost in the hydrodynamic limit. 

In~\cite{motsch2011numerical}, it was shown that the SOH model is the relaxation limit of a conservative system. Then, the authors proposed a splitting method by which the conservative part of the model was solved first, using a standard shock capturing scheme, and the relaxation part was solved next. The solutions obtained by this scheme provided very accurate approximations of the particle simulations~\cite{motsch2011numerical} or the kinetic ones~\cite{dimarcomotsch16}. It would be tempting to reproduce this strategy with the SOHB model and the first step is to find a conservative relaxation approximation of it. This is what we do in Section~\ref{sec_relaxation}. The application of this strategy to develop numerical simulations of the SOHB model will be deferred to future work. A special attention is devoted to the case of dimension 2 where the SOH and SOHB model coincide. However, we note that a SOH model with different coefficients (hence not coinciding with the  coefficients of the usual SOH model) can be derived from a higher-dimensional SOHB model under some invariance properties (a procedure which we call ``dimension reduction''). It results that different conservative relaxation approximations of this model can be constructed showing that such constructions are not unique in general. 

Another approach to deal with the ill-definition of discontinuous solutions of the SOHB model is to prevent them from occurring by introducing some viscosity. To this aim, we return to the kinetic description of the system. Indeed, macroscopic models are classically derived from kinetic equations when the ratio of the typical scales associated with the microscopic and macroscopic phenomena is sent to zero. For this, it is necessary to make scaling assumptions on the various dimensionless parameters of the model. In Section~\ref{sec:viscosity}, we show that, by a small modification of these scaling assumptions compared with~\cite{degond2021body}, we can keep viscous terms finite in the SOHB model. Our goal here is to provide the explicit form of these terms. The influence of this finite viscosity on the solutions will be studied in future work. We note that, even for the SOH model where viscosity terms were derived in~\cite{degond20133hydrodynamic}, the influence of the viscosity on the solution profiles has not been documented yet. Existence of solutions for the viscous SOH model was shown in~\cite{zhang2017local}. 

\smallskip
The main contributions of this work concern the SOHB model and are: 

(i) the proof of the hyperbolicity of the model, 

(ii) the derivation of a conservative relaxation approximation of the model, 

(iii) the derivation of viscous corrections to the model from the kinetic model. 

\smallskip
The organization of the paper is as follows. In Section~\ref{sec_SOHB_presentation}, we give a presentation of the model. Then, its hyperbolicity is shown in Section~\ref{sec_hyperbolic}. A conservative relaxation of the SOHB model is exhibited in Section~\ref{sec_relaxation}. Then, viscous corrections to the SOHB model are derived in Section~\ref{sec:viscosity}. Finally a conclusion and some perspectives are given in Section~\ref{sec:conclu}. Appendices~\ref{sec_invar_coor_change_proof} and~\ref{sec_dimension_reduc_proof} provide details of some technical proofs.

\setcounter{equation}{0}
\section{The Self-Organized Hydrodynamics model for body orientation}
\label{sec_SOHB_presentation}

The Self-Organized Hydrodynamics model for body orientation (SOHB) depicts the motion of a fluid consisting of a large number of small self-propelled rigid bodies in an~$n$-dimensional space. The attitude of a rigid-body can be described by a direct orthonormal frame attached to it (the agent's body frame), or equivalently by a rotation mapping a fixed direct orthonormal reference frame~$(e_1, \ldots, e_n)$ to the the agent's body frame. We will consider identical rigid-bodies with identical body frames attached to them. They are self-propelled with identical and uniform speed~$c_0$ and velocity directed along the first vector of their attached body frame. Since we are interested in a macroscopic description of this fluid, we are only interested in the average body frame in an infinitesimally small volume about an arbitrary position~$x \in {\mathbb R}^n$ at time~$t \geq 0$. How this average arises will be described in Section~\ref{sec:viscosity}, where the derivation of the model will be sketched. This average body frame is denoted by~$(\Omega_1(x,t), \ldots , \Omega_n(x,t))$. It is also a direct orthonormal frame and the rotation mapping the reference frame~$(e_1, \ldots, e_n)$ to it will be denoted by~$\rot(x,t)$ (i.e. 
$ \Omega_k = \rot e_k$,~$k = 1, \ldots, n$). Thus~$\rot(x,t)$ belongs to the group~$\textrm{SO}_n {\mathbb R}$ of~$n$-dimensional rotation matrices. In addition to giving the evolution of~$\rot$, the SOHB model gives the evolution of the mean density~$\rho = \rho(x,t)$ of agents around point~$x$ at time~$t$. 

The model is written as follows: 
\begin{align}
&
\partial_t \rho + \nabla \cdot (c_1 \rho \Omega_1) = 0,  \label{eq:mass_2} \\
&
\rho \big( \partial_t + c_2 \, \Omega_1 \cdot \nabla \big) \rot - {\mathbb W} \rot = 0,  
\label{eq:orient_3} 
\end{align}
where 
\begin{align}
\Omega_1 &= \rot e_1, \label{eq:def_Omega1} \\
{\mathbb W} &= F \wedge \Omega_1 - c_4 \, \rho \, \nabla \wedge \Omega_1,  \label{eq:expressW} \\
F &= - c_3 \, \nabla \rho - c_4 \, \rho \, \rot \nabla \cdot \rot. \label{eq:expressF} 
\end{align} 
The constants~$c_i$,~$i=1, \ldots, 4$ are such that~$c_1$,~$c_3 \in (0,\infty)$ and~$c_2$,~$c_4 \in {\mathbb R}$. As will be shown in Section~\ref{sec:viscosity}, this model can be derived from a kinetic model by means of a hydrodynamic limit. In such case, the coefficients~$c_i$ are linked with those of the kinetic model by explicit formulas (see e.g.~\cite{degondfrouvellemerino17, degondfrouvellemerinotrescases18, degond2018alignment} in dimension 3 and~\cite{degond2021body} in arbitrary dimension). However, here, our goal is to study the SOHB model~\eqref{eq:mass_2}-\eqref{eq:expressF} in full generality, without reference to a specific kinetic model (except in Section~\ref{sec:viscosity}), so the values of the coefficients will remain unspecified. The vector~$e_1$ is given (note indeed that the model is self-contained without reference to the other reference vectors~$e_2, \ldots, e_n$). 

Here, for two vector fields~$X = (X_i)_{i=1, \ldots, n}$ and~$Y = (Y_i)_{i=1, \ldots, n}$, 
$X \wedge Y$ and~$\nabla \wedge X$ stand for the antisymmetric matrices with entries 
\[ (X \wedge Y)_{ij} = X_i Y_j - X_j Y_i, \quad (\nabla \wedge X)_{ij} = \partial_{x_i} X_j - \partial_{x_j} X_i.\]
The quantity~$\nabla \cdot X$ is the divergence of~$X$:~$\nabla \cdot X = \sum_{i=1}^n \partial_{x_i} X_i$, while~$X \cdot \nabla$ stands for the operator~$X \cdot \nabla = \sum_{i=1}^n X_i \partial_{x_i}$. For a scalar field~$\varphi$,~$\nabla \varphi$ is the gradient vector whose components are~$(\nabla \varphi)_i = \partial_{x_i} \varphi$. For a type 2 tensor~$A = (A_{ij})_{i,j = 1 , \ldots, n}$ we denote by~$\nabla \cdot A$ its divergence, which is the vector of components~$(\nabla \cdot A)_j = \sum_{i=1}^n \partial_{x_i} A_{ij}$. The expression~$\rot \nabla \cdot \rot$ simply means the matrix~$\rot$ multiplying the vector~$\nabla \cdot \rot$. We have 
\begin{equation}
\rot \nabla \cdot \rot = \sum_{k=1}^n (\nabla \cdot \Omega_k) \Omega_k.
\label{eq:Thet_nab_Thet}
\end{equation} 

Eq.~\eqref{eq:mass_2} is nothing but the continuity equation for the fluid, and shows that the fluid velocity is given by~$c_1 \Omega_1$. Hence, the direction of the fluid velocity coincides with the first vector of the average body frame~$(\Omega_1, \ldots, \Omega_n)$ in the same way as the agent's self-propulsion velocity in the particle model is in the direction of the first vector of the agent's body frame. The fluid speed is equal to~$c_1$. 

The expression~$\partial_t + c_2 \, \Omega_1 \cdot \nabla$ at the left-hand side of~\eqref{eq:orient_3} is a convective derivative along the vector field~$c_2 \, \Omega_1$. Thus, it computes the time derivative of~$\Theta$ as one moves along the integral curves of this vector field. We note that, as~$c_2 \not = c_1$ in general,~$c_2 \, \Omega_1$ differs from the fluid velocity. We recall that for~$\rot \in \textrm{SO}_n {\mathbb R}$, the tangent manifold~$T_\rot$ to~$\textrm{SO}_n {\mathbb R}$ at~$\rot$ is given by~$T_\rot = \{ P \rot \, | \, P \in {\mathcal A} \}$, where~${\mathcal A}$ is the space of~$n \times n$ antisymmetric matrices. Since~$\big( \partial_t + c_2 \, \Omega_1 \cdot \nabla \big)$ is a derivation, the left-hand side of~\eqref{eq:orient_3} belongs to~$T_{\rot (x,t)}$ and thus is of the form~${\mathbb W} \rot$ with~${\mathbb W} \in {\mathcal A}$. We check that the expression~\eqref{eq:expressW} of~${\mathbb W}$ produces an antisymmetric matrix, showing the consistency of~\eqref{eq:orient_3}. In dimension~$n=3$ an antisymmetric matrix is realized as the cross product with a given vector and in rigid body dynamics, this vector is called the angular rotation vector. Thus, eq.~\eqref{eq:orient_3} is just rigid-body dynamics in dimension~$n$ and we can refer to~${\mathbb W}$ as the angular velocity matrix. 

Eq.~\eqref{eq:expressW} shows that~${\mathbb W}$ decomposes in two terms. The first one involves the vector~$F$ which is nothing but the force acting on the fluid. According to Eq.~\eqref{eq:expressF}, this force has two components. The first one is opposite to the density gradient and describes the effect of the pressure force. The second one is proportional to~$\nabla \cdot \rot$ and describes the effects of body-attitude gradients. While the first term is classical in fluid dynamics, the second one introduces effects specific to interacting rigid bodies. The antisymmetric matrix~$F \wedge \Omega_1$ describes an infinitesimal rotation in the plane spanned by~$\{\Omega_1,F\}$ (i.e. leaving the orthogonal of Span$\{\Omega_1,F\}$ invariant) which contributes to align~$\Omega_1$ with~$F$. Hence, due to the pressure component of~$F$, the fluid goes down the density gradients and due to its second component, it tends to escape the regions of large variations of body-attitude. In~\cite{degond2021body}, the second term of the expression~\eqref{eq:expressW} of~${\mathbb W}$ is shown to be written 
\[ -  c_4 \rho \nabla \wedge \Omega_1 = c_4 \rho \big( (\Omega_1 \cdot \nabla_x) \rot + {\mathbb A} \big),\]
where~${\mathbb A}$ is an antisymmetric matrix which leaves~$\Omega_1$ invariant (i.e.~${\mathbb A} \Omega_1 = 0$). Hence, this second term decomposes into transport along~$\Omega_1$ and rotation about the self-propulsion velocity~$\Omega_1$. We refer the reader to~\cite{degond2021body} for the expression of~${\mathbb A}$ which will not be used here. 

In the next section, we investigate the hyperbolicity and invariance properties of this system.

\setcounter{equation}{0}
\section{Hyperbolicity}
\label{sec_hyperbolic}

The goal of the present section is to prove the following

\begin{theorem}
The SOHB system~\eqref{eq:mass_2}-\eqref{eq:expressF} is weakly hyperbolic (i.e. it has real eigenvalues). In the generic case (i.e. except for special combinations of the coefficients and possibly a discrete set of propagation directions), it is hyperbolic (i.e. it is diagonalizable). In this generic case, it is strictly hyperbolic (i.e. it has distinct eigenvalues) in dimension~$n=3$ and not strictly hyperbolic in dimension~$n \geq 4$.
\label{eq:thm_hyperbolic}
\end{theorem}

The proof of this theorem is cumbersome unless one chooses appropriate references frames thanks to the model invariances. So, in Section~\ref{subsec_invariances} we state and prove these invariance properties, while the proof of Theorem~\ref{eq:thm_hyperbolic} is developed in Section~\ref{subsec_hyperbolic_proof}. 

\subsection{Invariances}
\label{subsec_invariances}

In this section, we show that the SOHB model enjoys two kinds of invariances. The first one is invariance under the choice of the reference basis~$(e_1, \ldots, e_n)$. The second one is invariance by coordinate change. We start with the first invariance

\begin{proposition}
The SOHB system~\eqref{eq:mass_2}-\eqref{eq:expressF} is invariant under a change of direct orthonormal reference frame~$(e_1, \ldots, e_n)$ (see the  proof below for the precise definition of this invariance). 
\label{prop:invariance_frame_change}
\end{proposition}

\textbf{Proof.}  Let~$(e'_1, \ldots, e'_n)$ be another direct orthonormal frame and~$R$ the matrix of the rotation which maps~$(e_1, \ldots, e_n)$ to~$(e'_1, \ldots, e'_n)$. Let~$\rot'$ be the matrix of the rotation which maps~$(e'_1, \ldots, e'_n)$ to~$(\Omega_1, \ldots, \Omega_n)$. Then, we obviously have~$\rot = \rot' R$ and~$\Omega_k = \rot e_k = \rot' e'_k$. One can also check that~$\Theta \nabla\cdot \Theta = \Theta'\nabla \cdot \Theta'$ using that~$R^TR=\textrm{I}$. Multiplying~\eqref{eq:orient_3} to the right by~$R$, we get that~$\rot'$ also satisfies~\eqref{eq:orient_3}. This shows that the system~\eqref{eq:mass_2}-\eqref{eq:expressF} is unchanged. \endproof

We now denote by~$(f_1, \ldots, f_n)$ the direct orthonormal frame associated to the coordinate system~$(x_1, \ldots, x_n)$, i.e. such that~$x=(x_1, \ldots, x_n)^T$ are the coordinates of the point~${\mathbf x} = \sum_{i=1}^n x_i \, f_i$ (the exponent~$T$ denotes the transpose of an array). We now have the  

\begin{proposition}
The SOHB system~\eqref{eq:mass_2}-\eqref{eq:expressF} is invariant under a change of direct orthonormal coordinate  frame~$(f_1, \ldots, f_n)$ (see Appendix~\ref{sec_invar_coor_change_proof} for the precise definition of this invariance).
\label{prop:invariance_coordinate_change}
\end{proposition}

\textbf{Proof.} The quantities involved in System~\eqref{eq:mass_2}-\eqref{eq:orient_3} are all scalars, vectors and type 2 tensors in the terminology of differential geometry. Using the classical formulas for coordinate changes of scalars, vectors and tensors immediately leads to the invariance of the system by change of direct orthonormal coordinate systems. However, for the unfamiliar reader to differential geometry, we give a direct proof below in Appendix~\ref{sec_invar_coor_change_proof}.

\subsection{Proof of hyperbolicity theorem~\ref{eq:thm_hyperbolic}}
\label{subsec_hyperbolic_proof}

\textbf{Proof.} Let~$(\rho_0,\rot_0) \in [0,\infty) \times \textrm{SO}_n{\mathbb R}$ be a uniform (i.e. independent of~$x$) state. Thanks to Prop.~\ref{prop:invariance_frame_change}, we may choose the reference frame~$(e_1, \ldots, e_n)$ such that~$\rot_0 = \textrm{I}$, where~$\textrm{I}$ is the identity matrix. To linearize~\eqref{eq:orient_3}, we write 
\[\rho = \rho_0 + \delta \sigma, \quad \rot = \exp ( \delta A), \quad \textrm{with} \quad \sigma =\sigma(x,t) \in {\mathbb R}, \, \, A=A(x,t) \in {\mathcal A} \quad \textrm{and} \quad \delta \ll 1,\]
which is legitimate as long as~$A$ stays in a bounded open set of~${\mathcal A}$ and~$\delta$ is small enough because the exponential map is a local diffeomorphism from a neighborhood of the origin in~${\mathcal A}$ to a neighborhood of~$\textrm{I}$ in~$\textrm{SO}_n{\mathbb R}$. We compute successively:  
\begin{align}
&
\rot = \textrm{I} + \delta A + {\mathcal O}(\delta^2), \nonumber \\
&
\Omega_k = \rot e_k = e_k + \delta \omega_k + {\mathcal O}(\delta^2), \textrm{ with } \omega_k = A e_k, \nonumber \\
&
\rot \nabla \cdot \rot = \sum_{k=1}^n (\nabla \cdot \Omega_k) \Omega_k = \delta \sum_{k=1}^n (\nabla \cdot \omega_k) e_k + {\mathcal O}(\delta^2), \label{eq:rotnarotlin} \\
&
F = \delta \Big( -c_3 \nabla \sigma - c_4 \rho_0 \sum_{k=1}^n (\nabla \cdot \omega_k) e_k \Big) + {\mathcal O}(\delta^2), \nonumber \\
&
{\mathbb W} = \delta \Big( -c_3 \nabla \sigma \wedge e_1  - c_4 \rho_0 \sum_{k=2}^n (\nabla \cdot \omega_k) e_k \wedge e_1 - c_4 \rho_0 \nabla \wedge \omega_1 \Big)  + {\mathcal O}(\delta^2), \nonumber \\
&
(\partial_t + c_2 \Omega_1 \cdot \nabla) \rot = \delta (\partial_t + c_2 e_1 \cdot \nabla) A + {\mathcal O}(\delta^2), \nonumber \\
&
\nabla \cdot (\rho \Omega_1) = \delta ( e_1 \cdot \nabla \sigma + \rho_0 \nabla \cdot \omega_1 ) + {\mathcal O}(\delta^2). \nonumber 
\end{align}
Thus, the pair~$(\sigma,A)$ satisfies the following linearized SOHB system: 
\begin{align}
&
\partial_t \sigma + c_1 (e_1 \cdot \nabla) \sigma + c_1 \rho_0 (\nabla \cdot \omega_1) = 0, \label{eq:mass_linearized} \\
&
\rho_0 (\partial_t A+ c_2 (e_1 \cdot \nabla) A)  = -c_3 \nabla \sigma \wedge e_1  - c_4 \rho_0 \sum_{k=2}^n (\nabla \cdot \omega_k) e_k \wedge e_1 - c_4 \rho_0 \nabla \wedge \omega_1, \label{eq:orient_linearized} 
\end{align}
with~$A \in {\mathcal A}$.

Now, we state the following invariance result for System~\eqref{eq:mass_linearized},~\eqref{eq:orient_linearized} and defer its proof to the end of this section. 

\begin{lemma}
System~\eqref{eq:mass_linearized},~\eqref{eq:orient_linearized} is invariant by change of~$(e_1, \ldots, e_n)$ to any other direct orthonormal reference frame~$(e'_1, \ldots, e'_n)$ such that~$e'_1 = e_1$. 
\label{lem:invar_linear}
\end{lemma}

We now take the Fourier transform with respect to~$x$ of~$\sigma$ and~$A$, defining~$\hat \sigma$,~$\hat A$ as functions of~$(\xi, t)$ where~$\xi \in {\mathbb R}^n$ is the Fourier dual variable of~$x$.
Also introducing~$\nu = \frac{\xi}{|\xi|}$ (using that~$\xi \not = 0$ because the case~$\xi=0$ is irrelevant for the definition of hyperbolicity, see below), we can write the Fourier transformed system as 
\begin{align}
&
\Big( \frac{1}{i |\xi|} \partial_t  + c_1 (e_1 \cdot \nu) \Big) \hat \sigma + c_1 \rho_0 (\nu \cdot \hat \omega_1) = 0, \label{eq:mass_Fourier} \\
&
\rho_0 \Big( \frac{1}{i |\xi|} \partial_t + c_2 (e_1 \cdot \nu) \Big) \hat A  = -c_3 (\nu \wedge e_1) \hat \sigma  - c_4 \rho_0 \sum_{k=2}^n (\nu \cdot \hat \omega_k) e_k \wedge e_1 - c_4 \rho_0 \nu \wedge \hat \omega_1. \label{eq:orient_Fourier} 
\end{align}

Using Lemma~\ref{lem:invar_linear}, we can rotate~$(e_2, \ldots, e_n)$ so as to bring~$e_2$ in the plane spanned by~$e_1$ and~$\nu$ (except if~$e_1$ and~$\nu$ are colinear, in which case we do not change~$(e_2, \ldots, e_n)$). Thus, we can write 
\[ \nu = \cos \theta e_1 + \sin \theta e_2,\]
with~$\theta \in [0,\pi)$ and we have 
\[e_1 \cdot \nu = \cos \theta, \qquad - \nu \wedge e_1 = e_1 \wedge \nu = \sin \theta \, e_1 \wedge e_2.\]
We can also write 
\begin{equation*} 
\hat A = \sum_{1 \leq k < \ell \leq n} {\hat A}_{k \ell} \, e_k \wedge e_\ell, 
\end{equation*}
since~$(e_k \wedge e_\ell)_{1 \leq k < \ell \leq n}$ is a basis of~$\mathcal{A}$
and we compute
\begin{align} 
\hat \omega_k &= \hat A e_k = \sum_{1 \leq \ell < m \leq n} \hat A_{\ell m} (e_\ell \wedge e_m) e_k 
= \sum_{1 \leq \ell < m \leq n} \hat A_{\ell m} (\delta_{mk} e_\ell - \delta_{\ell k} e_m ) \nonumber
\\
&= \sum_{\ell<k} \hat A_{\ell k} e_\ell - \sum_{m>k} \hat A_{km} e_m. \label{eq:hatomk}
\end{align}
Thus, we find 
\begin{align*}
\hat \omega_k \cdot \nu &= \cos \theta \Big( \sum_{\ell<k} \hat A_{\ell k} \delta_{\ell 1} - \sum_{m>k} \hat A_{km} \delta_{m1} \Big) 
+ \sin \theta \Big( \sum_{\ell<k} \hat A_{\ell k} \delta_{\ell 2} - \sum_{m>k} \hat A_{km} \delta_{m2} \Big) \\
&= \begin{cases}
- \hat A_{12} \sin \theta & \textrm{ if } k=1, \\
\hat A_{12} \cos \theta & \textrm{ if } k=2, \\
\hat A_{1k} \cos \theta + \hat A_{2k} \sin \theta & \textrm{ if } k \geq 3. 
\end{cases}
\end{align*}
Finally from~\eqref{eq:hatomk}, we get 
\[ \hat \omega_1 = - \sum_{m \geq 2} \hat A_{1m} e_m,\]
from which we deduce 
\begin{align*}
- \nu \wedge \hat \omega_1 &= \hat \omega_1 \wedge \nu =  - \Big(\sum_{m \geq 2} \hat A_{1m} e_m \Big) \wedge (\cos \theta e_1 + \sin \theta e_2) \\
&= \hat A_{12} \cos \theta e_1 \wedge e_2 +\sum_{m \geq 3} \hat A_{1m}  (\cos \theta \, e_1 \wedge e_m + \sin \theta \, e_2 \wedge e_m ). 
\end{align*}

Now, decomposing~\eqref{eq:orient_Fourier} on the basis~$(e_i \wedge e_j)_{1 \leq i < j \leq n}$, we can write System~\eqref{eq:mass_Fourier},~\eqref{eq:orient_Fourier} in the form 
\begin{equation}
 \frac{1}{i |\xi|} \partial_t U + {\mathbb M} U = 0, 
\label{eq:sys_gene}
\end{equation}
where~$U = (\hat \sigma, (\hat A_{ij})_{1 \leq i < j \leq n})^T$ is the vector of unknowns and~${\mathbb M}$ will be further described. We recall the following definitions: the SOHB System is: 
\begin{itemize}
    \item hyperbolic if and only if the matrix~${\mathbb M}$ is diagonalizable with real eigenvalues for all values of the unperturbed state~$(\rho_0, \rot_0)$ and all values of~$\xi \in {\mathbb R}^n$; 
    \item strictly hyperbolic if and only if it is hyperbolic and all eigenvalues are simple; 
    \item weakly hyperbolic if and only if all eigenvalues of~${\mathbb M}$ are real (but~${\mathbb M}$ is not necessarily diagonalizable). 
\end{itemize}
As we have seen, wlog, we can take~$\rot_0 = {\mathrm I}$, so it is enough that the hyperbolicity condition be true for all~$\rho_0 >0$. Likewise, the matrix~${\mathbb M}$ does not depend on~$|\xi|$ so it is enough that the hyperbolicity condition be true for all~$\nu \in {\mathbb S}^{n-1}$. We see that, if the SOHB system is hyperbolic, the solutions of the linearized system~\eqref{eq:sys_gene} have constant~$L^2$ norm in time. It is not true if the system is only weakly hyperbolic: the~$L^2$-norm may grow but only polynomially (by contrast to an unstable system in which the norm grows exponentially). The eigenvalues of~${\mathbb M}$ are called the characteristic speeds. They are the speeds of travelling-wave solutions of the linearized system corresponding to initial conditions belonging to the corresponding eigenspaces.

In fact, the system decomposes into the following uncoupled systems:
\begin{itemize}
\item A system for the pair~$(\hat \sigma, \hat A_{12})$ written as follows: 
\begin{align}
&
\frac{1}{i |\xi|} \partial_t \hat \sigma + c_1 \cos \theta \hat \sigma - c_1 \rho_0 \sin \theta \hat A_{12} = 0, \label{eq:sigA12_eq1} \\
&
\frac{1}{i |\xi|} \partial_t \hat A_{12} - \frac{c_3}{\rho_0} \sin \theta \hat \sigma + (c_2-2c_4) \cos \theta \hat A_{12}  =  0, \label{eq:sigA12_eq2}
\end{align}

\item Systems for the pairs~$(\hat A_{1k}, \hat A_{2k})$ with~$k$ such that~$3 \leq k \leq n$, written as follows: 
\begin{align}
&
\frac{1}{i |\xi|} \partial_t \hat A_{1k} + (c_2-2 c_4) \cos \theta \hat A_{1k} - c_4 \sin \theta \hat A_{2k} = 0
, \label{eq:A1kA2k_eq1} \\
&
\frac{1}{i |\xi|} \partial_t \hat A_{2k}  - c_4 \sin \theta \hat A_{1k} + c_2 \cos \theta \hat A_{2k} =0, \label{eq:A1kA2k_eq2}
\end{align}

\item Single equations for the quantities~$\hat A_{k\ell}$ for all~$(k,\ell)$ such that~$3 \leq k < \ell \leq n$, written as follows:  
\begin{equation}
\frac{1}{i |\xi|} \partial_t \hat A_{k\ell}  + c_2 \cos \theta \hat A_{k\ell} = 0. 
\label{eq:Akl}
\end{equation}
\end{itemize}

From this, we easily reconstruct the matrix~${\mathbb M}$ which is block diagonal, with one block of size~$2 \times 2$ corresponding to System~\eqref{eq:sigA12_eq1},~\eqref{eq:sigA12_eq2},~$n-2$ blocks of size~$2 \times 2$ corresponding to Systems~\eqref{eq:A1kA2k_eq1},~\eqref{eq:A1kA2k_eq2}, and~$(n-2)(n-3)/2$ blocks of size~$1 \times 1$ corresponding to Equations~\eqref{eq:Akl}. 

It is an easy matter to check that the block corresponding to System~\eqref{eq:sigA12_eq1},~\eqref{eq:sigA12_eq2} has eigenvalues 
\begin{equation*} \lambda_{\pm} = \frac{1}{2} \Big\{ (c_1 + c_2 - 2 c_4) \cos \theta \pm \Big( (c_1 - (c_2 - 2 c_4))^2 \cos^2 \theta + 4 c_1 c_3 \sin^2 \theta \Big)^{1/2} \Big\}. 
\end{equation*}
Since~$c_1 c_3 > 0$ by assumption, the quantity inside the square root is always nonnegative so that~$\lambda_\pm$ are real. Except if~$c_1 = c_2 - 2 c_4$ and~$\theta = 0$, these two eigenvalues are different and each of them is simple. 

Now, the blocks corresponding to System~\eqref{eq:A1kA2k_eq1},~\eqref{eq:A1kA2k_eq2} have eigenvalues 
\begin{equation*} 
\mu_{\pm} = (c_2 - c_4) \cos \theta \pm c_4, 
\end{equation*}
which again, are real. They are distinct as soon as~$c_4 \not = 0$. Each of these eigenvalues has multiplicity~$n-2$ as there are~$n-2$ such blocks. Finally, the single equation~\eqref{eq:Akl} has eigenvalue 
\begin{equation*}
\beta = c_2 \cos \theta,
\end{equation*}
which is again real. This eigenvalue has multiplicity equal to the number of such blocks, i.e.~$(n-2)(n-3)/2$. 

In the generic case~$\lambda_\pm$,~$\mu_\pm$ and~$\beta$ are pairwise distinct (there may be special combinations of the coefficients and values of the angle~$\theta$ for which two of them are equal). In this generic case, the matrix is diagonalizable. Indeed, identical eigenvalues correspond to distinct blocks, and each block is diagonalizable. Thus the dimension of the eigenspace is equal to the multiplicity of the corresponding eigenvalue. In this generic case, in dimension~$n=3$, there are no identical eigenvalues, and so, the system is strictly hyperbolic. However, in dimension~$n \geq 4$, there are identical eigenvalues and the system is not strictly hyperbolic. This ends the proof of Theorem~\ref{eq:thm_hyperbolic} \endproof

\textbf{Proof of Lemma~\ref{lem:invar_linear}.} The only term in System~\eqref{eq:mass_linearized},~\eqref{eq:orient_linearized} which depends on~$e_k$ for~$k=2, \ldots, n$ is the second term at the right-hand side of~\eqref{eq:orient_linearized}, i.e. the antisymmetric tensor 
\[{\mathbb T} = \tau \wedge e_1, \qquad \tau = \sum_{k=1}^n (\nabla \cdot \omega_k) e_k .\]
But from~\eqref{eq:rotnarotlin}, we also have~$\tau = \nabla \cdot A$ which does not depend on the choice of the reference frame~$(e_1, \ldots, e_n)$ and so,~$\tau$ is left invariant by any change of this frame. If the frame change leaves~$e_1$ untouched, then~${\mathbb T}$ itself is invariant, as are the other terms of System~\eqref{eq:mass_linearized},~\eqref{eq:orient_linearized}. This ends the proof. \endproof

\setcounter{equation}{0}
\section{Non-conservativity: first approach by conservative relaxation approximation}
\label{sec_relaxation}

\subsection{Rationale}
\label{subsec_relax_rationale}

The SOHB model is a hyperbolic model, but it is non-conservative. Indeed, while the continuity equation~\eqref{eq:mass_2} can be written in the form~$\partial_t \rho + \nabla \cdot \Phi_\rho = 0$ where~$\Phi_\rho = c_1 \rho \Omega_1$ is the mass flux, it is not obvious that Eq.~\eqref{eq:orient_3}  for~$\rot$ can be put in a similar form, i.e.~$\partial_t \rot_{ij} + \nabla \cdot \Phi_{ij} = 0$, for a convenient expression of~$\Phi_{ij}$. This situation was already encountered with the SOH model, which is the continuum version of the Vicsek alignment dynamics. The SOH system~\cite{degondfrouvelleliu15, degond20133hydrodynamic, degond2008continuum} describes the evolution of the density~$\rho(x,t)$ and mean self-propulsion direction~$\Omega (x,t)$ and is written
\begin{align}
&
\partial_t \rho + \nabla \cdot (c_1 \rho \Omega) = 0, \label{eq:mass_SOH} \\
&
\rho ( \partial_t + c_2 \Omega \cdot \nabla) \Omega + c_3 P_{\Omega^\bot} \nabla \rho = 0, \label{eq:orient_SOH}
\end{align}
where the constants~$c_1,\ c_2,\ c_3$ are generic constants (with~$c_1$,~$c_3 >0$ and~$c_2 \in {\mathbb R}$) which may be different from those of the SOHB model. The operator~$P_{\Omega^\perp}$ denotes the matrix of the orthogonal projection on the space orthogonal to~$\Omega$ and is writen~$P_{\Omega^\perp} = \textrm{I} - \Omega \otimes \Omega$. In the case of the SOH model, a method was empirically found to overcome the problem of the ill-definition of the model for discontinuous solutions (see discussion in the introduction). In~\cite{motsch2011numerical}, it was shown that the SOH model~\eqref{eq:mass_SOH},~\eqref{eq:orient_SOH} is the relaxation limit of a conservative model for the density~$\rho(x,t)$ and the fluid velocity~$v(x,t)$ given by
\begin{align}
&
\partial_t \rho + \nabla \cdot (\rho v) = 0, \label{eq:mass_relax} \\
&
\partial_t (\rho v) + \frac{c_2}{c_1} \nabla \cdot ( \rho v \otimes v) + c_1 c_3 \nabla \rho = \frac{1}{\alpha} \rho v (c_1^2 - |v|^2) ,  \label{eq:orient_relax}
\end{align}
where~$\alpha \ll 1$ is a relaxation parameter. Indeed, in the limit~$\alpha \to 0$, the relaxation term drives~$v$ towards a state where~$|v| = c_1$, i.e.~$v = c_1 \Omega$ where~$\Omega \in {\mathbb S}^{n-1}$, which leads to~\eqref{eq:mass_SOH}. Then, applying the operator~$P_{v^\bot}$ to~\eqref{eq:orient_relax} gets rid of the singular right-hand side and, taking the limit~$\alpha \to 0$, leads to~\eqref{eq:orient_SOH}. This was exploited numerically in~\cite{motsch2011numerical} to define a time-splitting scheme as follows. Supposing~$(\rho^k, \Omega^k) \approx (\rho(t^k), \Omega(t^k))$ is given. Then, advancing one time-step from~$t^k$ to~$t^k + \Delta t$ consists of:
\begin{itemize}
\item Step 1: solving the following system:   
\begin{align}
&
\partial_t \rho + \nabla \cdot (\rho v) = 0, \label{eq_mass_SOH_split}\\
&
\partial_t (\rho v) + \frac{c_2}{c_1} \nabla \cdot ( \rho v \otimes v) + c_1 c_3 \nabla \rho = 0 .  \label{eq_orient_SOH_split}
\end{align}
If we had~$c_2/c_1 = 1$, this system would coincide with the isothermal compressible Euler equations. One can easily extend standard shock-capturing schemes for the Euler equation~\cite{leveque1992numerical} to this system. This has been done in~\cite{motsch2011numerical} and later in~\cite{degond2013self, dimarcomotsch16}. In this step, the system is solved with initial datum~$(\rho^{k}, \Omega^{k})$ and leads to~$(\rho^{k+1}, v^{k+1}) = (\rho(\Delta t), v(\Delta t))$. 
\item Step 2: define~$\Omega^{k+1} = \frac{v^{k+1}}{|v^{k+1}|}$. 
\end{itemize}
The second step of the scheme corresponds to the limit~$\alpha \to 0$ of the solution to
\begin{align*}
&
\partial_t \rho  = 0,  \\
&
\partial_t (\rho v) = \frac{1}{\alpha} \rho v (c_1^2 - |v|^2),  
\end{align*}
which would be the natural second step of a standard splitting. Although no rigorous convergence result has been proved so far, this method has been observed to provide an extremely accurate match with the solutions of the particle system~\cite{dimarcomotsch16, motsch2011numerical}.

We note that System~\eqref{eq_mass_SOH_split},~\eqref{eq_orient_SOH_split} is hyperbolic if and only if~$c_2/c_1 \geq 1$~\cite{degond2013self}. However, this does not prevent the method to be used. Indeed, System~\eqref{eq_mass_SOH_split},~\eqref{eq_orient_SOH_split} only serves as an intermediate step in the overall splitting scheme used to approximate the hyperbolic System~\eqref{eq:mass_SOH},~\eqref{eq:orient_SOH}. Instabilities related to the non-hyperbolicity of this intermediate system only have a single time step to develop and then, are removed by the projection step. This feature has been exploited e.g. in~\cite{degond2013self} where the conservative system is solved by means of a shock capturing scheme which does not require the computation of the eigenvectors of the jacobian of the system (such as a polynomial scheme~\cite{cordier2014phase, degond1999polynomial}) and complex eigenvalues are replaced by their module. Simulations shown in~\cite{degond2013self} provide a validation of this approach. 

Now, a natural question is whether a similar methodology can be applied to the SOHB model and whether it leads to accurate approximations of the solutions of the underlying particle system. In the next section, we will partly answer this question positively by showing that the SOHB model can be obtained as the relaxation limit of a conservative problem in any dimension. The development of the corresponding splitting method and its assessment will be left for future work.

\subsection{Conservative relaxation approximation in arbitrary dimension }
\label{subsec_relax_dimn}

In this section, we show the following

\begin{proposition}
We assume that the reference frame~$(e_1, \ldots, e_n)$ coincides with the coordinate frame~$(f_1, \ldots, f_n)$. 
The following system, of unknowns~$(\rho, M)(x,t)$ where~$M$ is an~$n \times n$ matrix, is a relaxation approximation of the SOHB system~\eqref{eq:mass_2}-\eqref{eq:expressF}: 
\begin{align}
&
\partial_t \rho + c_1 \sum_{m=1}^n \partial_{x_m} (\rho M_{m1}) = 0, 
\label{eq:mass_relax_3} \\
&
\partial_t (\rho M_{ij}) - 2 c_4 \sum_{m=1}^n \partial_{x_m} (\rho M_{mj} M_{i1}) + c_2 \sum_{m=1}^n \partial_{x_m} (\rho M_{m1} M_{ij} ) \nonumber \\
& \hspace{8cm}
+ 2 (c_3-c_4) \partial_{x_i} \rho \, \delta_{j1} =  - \frac{1}{\alpha} \rho R_{ij} (M),
\label{eq:matrix_relax}
\end{align}
where the relaxation term~$R = (R_{ij})_{i,j = 1, \ldots, n}$ is given by 
\begin{equation}
R(M) = (MM^T - \textrm{I}) M.  
\label{eq:relax}
\end{equation}
In other words, if~$(\rho^\alpha, M^\alpha)$ is a solution to System~\eqref{eq:mass_relax_3}-\eqref{eq:relax} for a given~$\alpha >0$, then, formally, in the limit~$\alpha \to 0$,~$(\rho^\alpha, M^\alpha) \to (\rho,\rot)$ where~$\rot$ is an orthogonal matrix and~$(\rho, \rot)$ satisfies System~\eqref{eq:mass_2}-\eqref{eq:expressF}. Furthermore, if, for~$\alpha$ sufficiently small,~$\det M^\alpha$ stays positive, then~$\rot$ is a rotation. 
\label{prop:relax}
\end{proposition}

Defining~$v_i = M e_i$, System~\eqref{eq:mass_relax_3}-\eqref{eq:relax} is equivalent to the following system:
\begin{align}
&
\partial_t \rho + c_1 \nabla \cdot (\rho v_1) = 0, 
\label{eq:mass_relax_2} \\
&
\partial_t (\rho v_1) + (c_2 - 2 c_4) \nabla \cdot (\rho v_1 \otimes v_1) 
+ 2 (c_3-c_4) \nabla \rho  =  - \frac{1}{\alpha} \rho \big( \sum_{k=1}^n v_k \otimes v_k - \textrm{I} \big) v_1.
\label{eq:v1} \\
&
\partial_t (\rho v_j) - 2 c_4 \nabla \cdot (\rho v_j \otimes v_1) + c_2 \nabla \cdot (\rho v_1 \otimes v_j) \nonumber \\
& \hspace{6cm}
= - \frac{1}{\alpha} \rho \big( \sum_{k=1}^n v_k \otimes v_k - \textrm{I} \big) v_j, \quad j = 2, \ldots, n.
\label{eq:vj}
\end{align}
In this way, the relaxation model appears as an~$n$-velocity compressible fluid dynamics model. All quantities showing up in the model are scalars, vectors or tensors in the sense of differential geometry, i.e. they are changed through a change of coordinate frame~$(f_1, \ldots, f_n)$ to~$(f'_1, \ldots, f'_n)$ like scalars, vectors or tensors. It follows that the system~\eqref{eq:mass_relax_2}-\eqref{eq:vj} is invariant by change of the  coordinate frame~$(f_1, \ldots, f_n)$ (or equivalently, one may develop the same arguments as those of Appendix~\ref{sec_invar_coor_change_proof} and find that the system is written in the same way with any choice of coordinate frame~$(f_1, \ldots, f_n)$). Thus, in~\eqref{eq:mass_relax_2}-\eqref{eq:vj}, the frames~$(e_1, \ldots, e_n)$ and~$(f_1, \ldots, f_n)$ are not constrained to coincide. The matrix~$M$ can also be taken as the matrix which maps~$(e_1, \ldots, e_n)$ to~$(v_1, \ldots, v_n)$ in any coordinate frame, but then~\eqref{eq:matrix_relax} must be transformed accordingly as the relation~$(Me_j)_i = M_{ij}$ is not true anymore. In the same way as for the SOHB system, the resulting equation will be invariant by changing~$M$ to~$M' = MR$ where~$R$ is a rotation (or more generally an orthogonal matrix). We will not use the equation for~$M$ in an arbitrary coordinate frame, as we will prefer to write the system in the form~\eqref{eq:mass_relax_2}-\eqref{eq:vj}. 

It is interesting to note that, in the SOHB system~\eqref{eq:mass_2}-\eqref{eq:expressF}, the same coefficient~$c_4$ appears in two different places (multiplying the term~$\rho \rot \nabla \cdot \rot \wedge \Omega_1$ and multiplying~$\rho \nabla \wedge \Omega_1$). The system would make sense if these two terms were multiplied by a different coefficient. However, in the relaxation system~\eqref{eq:mass_relax_3}-\eqref{eq:relax},~$c_4$ only appears alone in a single place, namely multiplying the second term of~\eqref{eq:matrix_relax} (in the fourth term,~$c_4$ is combined with~$c_3$ which is arbitrary anyway). So, if the SOHB system had different coefficients in front of the two term involving~$c_4$, it could not be obtained as the relaxation limit of a system of the type~\eqref{eq:mass_relax_3}-\eqref{eq:relax}. It looks like the SOHB system has the exact structure needed to be the relaxation limit of a conservative system. This may be not coincidental but rather, the sign that the relaxation system~\eqref{eq:mass_relax_3}-\eqref{eq:relax} is itself the hydrodynamic limit of a suitable particle system. This point will be investigated in future work, as well as the development of numerical approximations of the SOHB system based on this relaxation model using the strategy of Section~\ref{subsec_relax_rationale}. 

\medskip
\noindent
\textbf{Proof.} 
We first study the homogeneous system 
\begin{equation} 
\frac{d M}{d t}  = - \frac{1}{\alpha} R (M), \qquad M(0) = M_0, 
\label{eq:relax_syst}
\end{equation}
and show that, if~$\det M_0 \not = 0$ and~$t>0$, 
\begin{equation*}
M(t) \to \rot_0 \quad \textrm{ as } \quad \alpha \to 0, \quad \textrm{ where } \quad \rot_0 \in \textrm{O}_n{\mathbb R}. 
\end{equation*} 
Furthermore,~$\rot_0$ is the polar decomposition of~$M_0$, i.e.~$M_0$ is uniquely factorized as~$M_0 = S_0 \rot_0$ where~$S_0$ is a symmetric positive definite matrix. We recall that~$S_0 = Q_0^{1/2}$ with~$Q_0 = M_0 M_0^T$ and~$\rot_0 = S_0^{-1} M_0$ where~$Q_0^{1/2}$ is defined as follows: since~$Q_0$ is symmetric and positive definite, there exists~$U \in \textrm{SO}_n{\mathbb R}$ and a diagonal matrix~$D_0$ such that~$Q_0 = U D_0 U^T$. Let~$d_{10}, \ldots, d_{n0}$ be the diagonal elements of~$D_0$. Then~$d_{i0} >0$,~$\forall i \in \{1, \ldots, n\}$ and~$Q_0^{1/2} = U D_0^{1/2} U^T$ where~$D_0^{1/2}$ is the diagonal matrix with diagonal elements~$\sqrt{d_{10}}, \ldots, \sqrt{d_{n0}}$. As a consequence, if~$\det M_0 > 0$,~$\rot_0  \in \textrm{SO}_n{\mathbb R}$. 

We first show that~$Q(t) =: M M^T(t) \to \textrm{I}$ as~$\alpha \to 0$. Indeed,~$Q$ satisfies the following system
\[ \frac{d Q}{d t}  = - \frac{2}{\alpha} (Q^2 - Q), \qquad Q(0) = Q_0.\]

Let~$D(t)$ be the solution of 
\[ \frac{d D}{d t}  = - \frac{2}{\alpha} (D^2 - D), \qquad D(0) = D_0,\]
and define~$\tilde Q(t) = U D(t) U^T$. We check that~$\tilde Q$ satisfies the same ordinary differential equation as~$Q$ with the same initial condition so, by the uniqueness statement of the Cauchy-Lipschitz theorem,~$\tilde Q =Q$. Now let~$d_1(t), \ldots, d_n(t)$ be the diagonal elements of~$D(t)$. They satisfy 
\begin{equation}
 \frac{d}{d t} d_i = - \frac{2}{\alpha} (d_i^2 - d_i), \qquad d_i(0) = d_{i0}. 
\label{eq:di}
\end{equation}
The solution of~\eqref{eq:di} is written
\[ d_i(t) = \frac{d_{i0}}{e^{-2t/\alpha} + d_{i0} (1 - e^{-2t/\alpha})}.\]
Now, when~$\alpha \to 0$ and~$t >0$,~$d_i(t) \to 1$ because~$d_{i0} > 0$. Thus,~$Q(t) \to \textrm{I}$.

Then, we consider~$S(t) = Q(t)^{1/2} = U D(t)^{1/2} U^T$ where~$D(t)^{1/2}$ is defined in a similar way as~$D_0^{1/2}$. Then,~$S^2=Q$ and differentiating, we get 
\[ 2 \frac{d S}{d t} S = - \frac{2}{\alpha} (MM^T - \textrm{I}) S^2,\]
where we have used that~$S$ and~$\frac{dS}{dt} = U \frac{d}{dt} (D(t)^{1/2}) U^T$ commute as~$D(t)^{1/2}$ and~$\frac{d}{dt} (D(t)^{1/2})$ commute. Multiplying this equation to the right by~$S^{-1} S_0^{-1} M_0$, we end up with 
\[ \frac{d}{d t} (S S_0^{-1} M_0) = - \frac{1}{\alpha} (MM^T - \textrm{I}) (S S_0^{-1} M_0),  \qquad (S S_0^{-1} M_0)|_{t=0} = M_0.\]
Thus,~$S S_0^{-1} M_0$ and~$M$ satisfy the same equation with the same initial condition. Again, by uniqueness, they are equal. Hence~$M = S S_0^{-1} M_0 = S \rot_0$. Now, since~$Q \to \textrm{I}$ as~$\alpha \to 0$, we have~$S \to \textrm{I}$ and thus,~$M \to \rot_0$, which ends the proof of this step. 

Now, we go back to System~\eqref{eq:mass_relax_3}-\eqref{eq:relax} and denote by~$(\rho^\alpha, M^\alpha)$ its solution for a given value of~$\alpha$. Thanks to the properties of~\eqref{eq:relax_syst}, and assuming that~$\det M^\alpha(x,t) >0$ for all~$(x,t) \in {\mathbb R}^n \times [0,\infty)$ and all~$\alpha$ small enough, we formally have~$(\rho^\alpha, M^\alpha) \to (\rho, \rot) \in [0,\infty) \times \textrm{SO}_n{\mathbb R}$ as~$\alpha \to 0$.  Indeed, stable stationary states of the homogeneous system~\eqref{eq:relax_syst} corresponding to initial data with positive determinants belong to~$\textrm{SO}_n{\mathbb R}$. Thus, in the spatially non-homogenous case, solutions with~$\det M^\alpha >0$ will generically converge to~$\rot \in \textrm{SO}_n{\mathbb R}$ when~$\alpha \to 0$ because, at least formally, they are small perturbations of order~$\alpha$ of the spatially homogeneous system. We define~$\Omega_i = \rot e_i$. Then, the limit of~\eqref{eq:mass_relax_3} is clearly~\eqref{eq:mass_2}. The difficulty is in showing that~$\rot$ satisfies~\eqref{eq:orient_3} with~${\mathbb W}$ given by~\eqref{eq:expressW}. We write~\eqref{eq:matrix_relax} as 
\[ {\mathcal T}(\rho^\alpha, M^\alpha) = - \frac{1}{\alpha} \rho^\alpha R(M^\alpha),\]
where~${\mathcal T}(\rho^\alpha, M^\alpha)$ is a short-hand notation for the left-hand side of~\eqref{eq:matrix_relax}. We denote by~${\mathcal M}$,~${\mathcal A}$,~${\mathcal S}$ the spaces of~$n \times n$ matrices over~${\mathbb R}$, of antisymmetric matrices and of symmetric matrices respectively. We also endow~${\mathcal M}$ with the euclidean structure generated by the Frobenius inner product~$A \cdot B = \textrm{Tr} \{ A^T B \}$. With this structure,~${\mathcal A}$ and~${\mathcal S}$ are orthogonal supplements of each other in~${\mathcal M}$. Then,~$R(M) \in {\mathcal S} \, M$. Thus, 
\[  P_{({\mathcal S} \, M^\alpha)^\bot} \big( {\mathcal T}(\rho^\alpha, M^\alpha) \big) = - \frac{1}{\alpha} \rho^\alpha P_{({\mathcal S} \, M^\alpha)^\bot}(R(M^\alpha)) = 0,\]
where~$P_{({\mathcal S} \, M^\alpha)^\bot}$ denotes the orthogonal projection from~${\mathcal M}$ onto~$({\mathcal S} \, M^\alpha)^\bot$. Taking the limit~$\alpha \to 0$ leads to 
\[  P_{({\mathcal S} \, \rot)^\bot} \big( {\mathcal T}(\rho, \rot) \big) = 0.\]
Now, it is readily seen that 
$({\mathcal S} \, \rot)^\bot = {\mathcal A} \, \rot = T_\rot$ where~$T_\rot$ is the tangent space to~$\textrm{SO}_n{\mathbb R}$ at~$\rot$, and that for any~$X \in {\mathcal M}$, 
\begin{equation}
P_{T_\rot} X = P_{{\mathcal A} \, \rot} X = \frac{X \rot^T - \rot X^T}{2} \rot. 
\label{eq:def_Ptheta}
\end{equation}
Thus, it is now sufficient to show that 
\begin{align}
\frac{1}{2} ( {\mathcal T}(\rho, \rot) \rot^T - \rot {\mathcal T}(\rho, \rot)^T ) = \rho &\partial_t \rot \, \rot^T  + c_2 \rho (\Omega_1 \cdot \nabla) \rot \, \rot^T \nonumber \\
& + c_3 \nabla \rho \wedge \Omega_1 + c_4 \rho \rot \nabla \cdot \rot \wedge \Omega_1 + c_4 \rho \nabla \wedge \Omega_1. \label{eq:match}
\end{align}
We label~${\mathcal T}_1$, \ldots,~${\mathcal T}_4$ the four terms at the left-hand side of~\eqref{eq:matrix_relax} (with~$M$ replaced by~$\rot$) and~${\mathcal U}_1$, \ldots,~${\mathcal U}_5$ the five terms at the right-hand side of~\eqref{eq:match}, in the order they appear in both cases. We 
also define 
\[ \tilde {\mathcal T}_k = \frac{1}{2} ( {\mathcal T}_k \rot^T - \rot {\mathcal T}_k^T ), \qquad k=1, \ldots, 4. \]
We compute successively
\[
\tilde {\mathcal T}_1 = \frac{1}{2} \big( \partial_t (\rho \rot) \, \rot^T - \rot \, \partial_t (\rho \rot)^T \big) = \frac{\rho}{2} \, ( \partial_t \rot \, \rot^T - \rot \, \partial_t \rot^T ) = \rho \, \partial_t \rot \, \rot^T = {\mathcal U}_1, 
\]
because~$ \partial_t (\rot \rot^T) = \partial_t \textrm{I} = 0 = \partial_t \rot \, \rot^T + \rot \, \partial_t \rot^T$. Then, 
\begin{align*}
(\tilde {\mathcal T}_2)_{ij} &= - c_4 \sum_{k,m = 1}^n \big( \partial_{x_m} (\rho \, \rot_{mk} \rot_{i1}) \rot_{jk} - \partial_{x_m} (\rho \, \rot_{mk} \rot_{j1}) \rot_{ik} \big) \\
&= -c_4 \sum_{k,m = 1}^n \begin{split}& \\ \big\{ \partial_{x_m} \rho \, \rot_{mk} \big(  \rot_{i1} \rot_{jk} - \rot_{j1} \rot_{ik} \big)
+ \rho \,  \partial_{x_m} \rot_{mk} \big( \rot_{i1} \rot_{jk} - \rot_{j1} \rot_{ik} \big)& \\
+ \rho \, \rot_{mk} \big( \partial_{x_m} \rot_{i1} \rot_{jk} - \partial_{x_m}\rot_{j1} \rot_{ik} &\big) \big\}\end{split}\\
& \mbox{} \\
&= -c_4 ( \partial_{x_j} \rho \, \rot_{i1}  - \partial_{x_i} \rho \, \rot_{j1} )
-c_4 \rho \big( \rot_{i1} (\rot \nabla \cdot \rot)_j - \rot_{j1} (\rot \nabla \cdot \rot)_i \big)) \\
& \hspace{8cm}
-c_4 \rho ( \partial_{x_j} \rot_{i1}  - \partial_{x_i} \rot_{j1} ) \\
& \mbox{} \\
&= c_4 (\nabla \rho \wedge \Omega_1)_{ij} + c_4 \rho \, (\rot \nabla \cdot \rot \wedge \Omega_1)_{ij} + c_4 \rho \, (\nabla \wedge \Omega_1)_{ij}. 
\end{align*}
Here, we have used that~$\sum_{k=1}^n \rot_{mk} \rot_{jk} = \delta_{mj}$ (and the same replacing~$j$ by~$i$) and that~$\rot_{i1} = (\Omega_1)_i$ (because here, we assume that the frames~$(e_1, \ldots, e_n)$ and~$(f_1, \ldots, f_n)$ coincide). 
Hence~$ \tilde {\mathcal T}_2 = c_4 \nabla \rho \wedge \Omega_1 + {\mathcal U}_4 + {\mathcal U}_5$. Now, we turn to 
\begin{align*}
(\tilde {\mathcal T}_3)_{ij} &= \frac{c_2}{2} \sum_{k,m = 1}^n \big( \partial_{x_m} (\rho \, \rot_{m1} \rot_{ik}) \rot_{jk} - \partial_{x_m} (\rho \, \rot_{m1} \rot_{jk}) \rot_{ik} \big) \\
&= \frac{c_2}{2}  \rho \,  \sum_{k,m = 1}^n  \rot_{m1} \big(  \partial_{x_m} \rot_{ik} \rot_{jk} - \partial_{x_m} \rot_{jk} \rot_{ik} \big) \\
&= c_2 \rho \, \sum_{k,m = 1}^n  \rot_{m1} \rot_{jk} \partial_{x_m} \rot_{ik} = c_2 \rho \big( (\Omega_1 \cdot \nabla) \rot \, \rot^T \big)_{ij},  
\end{align*}

where we have used that 
\[ \sum_{k=1}^n (\partial_{x_m} \rot_{ik} \rot_{jk} + \partial_{x_m} \rot_{jk} \rot_{ik} ) = \partial_{x_m} (\sum_{k=1}^n \rot_{ik} \rot_{jk}) = \partial_{x_m} \delta_{ij} = 0.\]
Hence,~$\tilde {\mathcal T}_3 = {\mathcal U}_2$. The final term is 
\begin{align*}
(\tilde {\mathcal T}_4)_{ij}
&= (c_3-c_4)  \sum_{k= 1}^n \big( \partial_{x_i} \rho \, \delta_{k1} \, \rot_{jk} - \partial_{x_j} \rho \, \delta_{k1} \, \rot_{ik} \big) \\
&= (c_3-c_4) \big( \partial_{x_i} \rho \, \rot_{j1} - \partial_{x_j} \rho \,  \rot_{i1} \big) = (c_3-c_4) (\nabla \rho \wedge \Omega_1)_{ij}, 
\end{align*}
so,~$\tilde {\mathcal T}_4 = {\mathcal U}_3 - c_4 \nabla \rho \wedge \Omega_1$. Summing all the~$\tilde {\mathcal T}_i$, we find that~\eqref{eq:match} is true, which ends the proof. \endproof

\subsection{Conservative relaxation approximation in dimension 2}
\label{subsec_relax_dim2}

Intuitively, in dimension~$n=2$, the SOHB and SOH model should coincide. Indeed, going back to the microscopic picture (see e.g.~\cite{degond2021body}) one would anticipate that aligning one vector of a two dimensional frame with the neighboring corresponding vectors, or aligning the two vectors of this frame with the neighboring frames should lead to the same result. However, previous works have only focused on the case~$n=3$~\cite{degondfrouvellemerino17, degondfrouvellemerinotrescases18, degond2018alignment} or~$n \geq 3$~\cite{degond2021body}. The reason why~$n=2$ was excluded from~\cite{degond2021body} is that~\cite{degond2021body} relies on representation theory of semi-simple groups and~$\textrm{SO}_2{\mathbb R}$ is abelian hence not semi-simple. However, it is easy to investigate the main  lemma depending on representation theory in~\cite{degond2021body} (Lemma 4.9) and, by direct computation, check that part (i) of the lemma still holds in dimension 2 while part (ii) also holds but the expression in factor of~$C_4$ vanishes identically. It follows that the SOHB model in dimension 2 is again given by~\eqref{eq:mass_2},~\eqref{eq:orient_3}, but with~$c_4 = 0$ in the expression of~${\mathbb W}$. Then, due to the dimension~$n=2$, ~$\Omega_2$ is just obtained from~$\Omega_1$ by a rotation of angle~$\pi/2$. So, letting~$\Omega = \Omega_1$, we can reduce the SOHB system to a system for~$(\rho, \Omega)$. After some easy computations, this system turns out to coincide with the SOH system~\eqref{eq:mass_SOH},~\eqref{eq:orient_SOH} as anticipated. 

However, there is another way to derive a SOHB system in dimension~$n=2$, which consists of starting from the SOHB system in higher dimension and supposing some uniformity of~$\rho$ and of the frame~$(\Omega_1, \ldots, \Omega_n)$ in a subspace of lower dimension. This process will be called dimension reduction. First, in Section~\ref{subsubsec_dimensional reduction}, we show how dimension reduction works in general. Then, in Section~\ref{subsubsec_dimensional reduction_n=2}, we apply it to derive a SOHB model in dimension~$n=2$ from the SOHB model in higher dimension. We will see that this ``dimensionly-reduced SOHB model in dimension 2'' is  the same as the SOH model obtained by the direct derivation discussed above, but with different coefficients.

\subsubsection{Dimension reduction in arbitrary dimension}
\label{subsubsec_dimensional reduction}

Here, we investigate under which conditions the SOHB system in dimension~$n$ reduces to an SOHB system in dimension~$p$ with~$2 \leq p \leq n-1$. As before, we are given a direct orthonormal reference frame~$(e_1, \ldots, e_n)$ for the definitions of the rotations. We denote by~$(\rho^0, \rot^0)$ the initial condition to the SOHB system~\eqref{eq:mass_2}-\eqref{eq:expressF} and assume existence and uniqueness of a smooth solution associated with this initial datum on a time interval~$[0,T]$. We also assume that~$\rho >0$ on this time interval. We denote by~$(\rho,\rot)$ this solution. We let~$\Omega_k^0 = \rot^0(e_k)$,~$k = 1, \ldots, n$. 

We let~$(f_1, \ldots, f_n)$ be a direct orthonormal coordinate system of~${\mathbb R}^n$ and the point~$x= \sum_{i=1}^n x_i f_i$ is identified with the~$n$-tuple of coordinates~$(x_1, \ldots, x_n)$. We will denote by~$\Pi$ the projection~${\mathbb R}^n \to {\mathbb R}^p$,~$(x_1, \ldots, x_n) \mapsto (x_1, \ldots x_p)$ and by~$i$ the canonical injection~${\mathbb R}^p \to {\mathbb R}^n$,~$\bar x =(x_1, \ldots x_p) \mapsto i(\bar x) = (x_1, \ldots, x_p, 0, \ldots, 0)$. 

Now, we assume that the following holds:

\begin{hypothesis}~

\begin{itemize}
\item[(i)]~$\Omega_{p+1}^0$,~$\Omega_{p+2}^0$, \ldots,~$\Omega_n^0$ are independent of~$x$. 
\item[(ii)] We have~$ \textrm{Span} \{ \Omega_{p+1}^0, \Omega_{p+2}^0, \ldots, \Omega_n^0 \} = \textrm{Span} \{ f_{p+1}, f_{p+2}, \ldots, f_n \}$.
\item[(iii)]~$\rho^0$ and~$\rot^0$ only depend on~$\Pi x$, i.e. there exists~$(\bar \rho^0,\bar \rot^0)$,~${\mathbb R}^p \to [0,\infty) \times \textrm{SO}_n{\mathbb R}$ such that~$ \rho^0(x) = \bar \rho^0(\Pi x), \, \rot^0(x) = \bar \rot^0(\Pi x)$.
\end{itemize}
\label{hyp:dim_reduc}
\end{hypothesis}
Note that we have chosen the last~$n-p$ indices for convenience in (i)-(iii) but the result would apply to any subset of~$n-p$ indices in the set~$\{2, \ldots, n\}$. Now, we have the

\begin{proposition} Under the assumption of local existence in time of a smooth solution for the SOHB system~\eqref{eq:mass_2}-\eqref{eq:expressF} in dimension~$p$ (for any smooth initial condition), we have the following properties :

(i) The solution~$(\rho, \rot)$ does not depend on the variables~$(x_{p+1}, \ldots, x_n)$, i.e. there exists~$(\bar \rho, \bar \rot):{\mathbb R}^p \times [0,T] \to [0,\infty) \times \textrm{SO}_n{\mathbb R}$ such that~$\rho(x,t) = \bar \rho(\Pi x,t), \, \rot(x,t) = \bar \rot(\Pi x,t)$. Furthermore~$\Omega_j$ is independent of time and space, equal to~$\Omega_j^0$ for all~$j>p$. 

(ii) Let~$(e'_1, \ldots, e'_n)$ be a new direct orthonormal reference frame for the definitions of the rotations, such that~$\textrm{Span} \{ e'_1, \ldots, e'_p \} = \textrm{Span} \{ f_1, \ldots, f_p \}$. Let~$R$ be the rotation that takes~$(e'_1, \ldots, e'_n)$ to~$(e_1, \ldots, e_n)$, i.e. such that~$e_j = R e'_j$ (for~$1\leq j\leq n$). Define~$\bar e'_j = \Pi e'_j$ and note that $e'_j = i (\bar e'_j)$. Let~$\theta$:~${\mathbb R}^p \times [0,T] \to {\mathcal M}_p {\mathbb R}$ (where 
${\mathcal M}_p {\mathbb R}$ is the space of real~$p \times p$ matrices), defined by~$\theta = \Pi \bar \rot R i$, i.e.~$\theta \bar e_j'(\Pi x,t)=\bar \rot(\Pi x,t) e_j=\Omega_j(x,t)$ for~$1\leq j\leq p$. Then,~$\theta \in \textrm{SO}_p{\mathbb R}$ and~$(\bar \rho, \theta)$ is a solution of the SOHB system~\eqref{eq:mass_2}-\eqref{eq:expressF} in dimension~$p$, associated with the orthonormal reference frame~$(\bar e'_1, \ldots, \bar e'_p)$ of~${\mathbb R}^p$. 

\label{prop_dim_reduc}
\end{proposition}

Note that the introduction of a new reference frame for the definition of the rotations $(\bar e'_1, \ldots, \bar e'_p)$ is needed because the original frame  $(e_1, \ldots, e_n)$ may not satisfy $\textrm{Span} \{ e_1, \ldots, e_p \} = \textrm{Span} \{ f_1, \ldots, f_p \}$. The rationale for this proposition comes from the fact that the conditions in Hypothesis~\ref{hyp:dim_reduc} are formally propagated in time by the SOHB system~\eqref{eq:mass_2}-\eqref{eq:expressF}, with the unit vectors~$\Omega_{p+1}$,~$\Omega_{p+2}$, \ldots,~$\Omega_n$ being actually constant in time. To see this we first write system~\eqref{eq:mass_2}-\eqref{eq:expressF} in terms of~$\rho$ and~$\Omega_j$,~$j=1, \ldots, n$:
\begin{align}
&
\partial_t \rho + \nabla \cdot (c_1 \rho \Omega_1) = 0, \label{eq:mass_rhom} \\
&
\rho ( \partial_t + (c_2-c_4) \Omega_1 \cdot \nabla) \Omega_1 = - c_3 P_{\Omega_1^\bot} \nabla \rho - c_4 \rho \sum_{k=2}^n (\nabla \cdot \Omega_k) \Omega_k , \label{eq:om1} \\
&
\rho ( \partial_t + c_2 \Omega_1 \cdot \nabla) \Omega_j = \Big( c_3 \Omega_j \cdot \nabla \rho + c_4 \rho \nabla \cdot \Omega_j \Big) \Omega_1 \nonumber \\
&\hspace{4cm}
+ c_4 \rho \Big( (\nabla \Omega_j) \Omega_1 + (\Omega_j \cdot \nabla) \Omega_1 \Big), \quad \forall j \in \{2, \ldots, n \}, \label{eq:omj}
\end{align}
where, for a vector field~$X$, we write~$(\nabla X)_{ij}=\partial_{x_i}X_j$. To get these equations, we have used that~$\Omega_1\cdot \Omega_j=0$, together with the following identities, for vector fields~$X$,~$Y$, and~$Z$:
\begin{align*}
 & (X\wedge Y)Z=(Y\cdot Z)X-(X\cdot Z)Y,\\
 &(\nabla\wedge X)Y=(\nabla X)Y-(Y\cdot\nabla)X=\nabla(X\cdot Y)-(\nabla Y)X-(Y\cdot\nabla)X,
\end{align*}
this last identity implying that~$(\nabla\wedge \Omega_1)\Omega_1=-(\Omega_1\cdot\nabla)\Omega_1$.

 Now, for~$j > p$, if the conditions of Hypothesis~\ref{hyp:dim_reduc} are valid until time~$t$, by the condition~$(i)$,~$\Omega_j$ is independent of~$x$ and thus,~$(\Omega_1 \cdot \nabla) \Omega_j = 0$,~$\nabla \cdot \Omega_j = 0$ and~$\nabla \Omega_j = 0$. We also have~$(\Omega_j \cdot \nabla) \Omega_1 =-(\Omega_1 \cdot \nabla) \Omega_j= 0$.
 Now, by the conditions~$(ii)$-$(iii)$, we obtain~$\Omega_j \cdot \nabla \rho = 0$ because~$\nabla \rho$ is contained in~$\textrm{Span} \{ f_1, \ldots, f_p \}$ and~$\Omega_j$ is orthogonal to it. Hence~$\rho \partial_t \Omega_j = 0$, which implies~$\partial_t \Omega_j = 0$ as~$\rho >0$ by assumption. Then, it is clear that~$\partial_t \rho$ and~$\partial_t \Omega_j$, for~$1 \leq j \leq p$ only depend on~$\Pi x$ as the right-hand sides of~\eqref{eq:mass_rhom}-\eqref{eq:omj} involve quantities that are only depending on~$\Pi x$. Thus, Assumptions~\ref{hyp:dim_reduc} are formally propagated in time, which corresponds to the first part of Proposition~\ref{prop_dim_reduc}, allowing to define~$\bar \rho$ and~$\bar \rot$ as functions of~$\bar x \in {\mathbb R}^p$ and time~$t$, the functions~$\bar{\Omega}_j=\bar{\rot}e_j$ for~$1\leq j\leq p$ satisfying equations~\eqref{eq:mass_rhom}-\eqref{eq:omj}, which translated in the frame~$(\bar{e}_1',\ldots,\bar {e}_p')$ corresponds to the second part of Proposition~\ref{prop_dim_reduc}. The rigorous way to prove Proposition~\ref{prop_dim_reduc} consists in constructing a solution of the reduced SOHB system in dimension~$p$ and extend it in dimension~$n$ by setting~$\Omega_j$ constant in time and space for~$j\geq p+1$, proving that this is a solution of the SOHB system in dimension~$n$ (and therefore the solution corresponding to the given initial conditions, by uniqueness). This construction is deferred to Appendix~\ref{sec_dimension_reduc_proof}. 
 We note that in the dimension reduction process, the coefficients~$c_1$, \ldots,~$c_4$ of the reduced dimension system are the same as those of the initial system.

\subsubsection{Dimension reduction to dimension 2 and conservative relaxation approximation}
\label{subsubsec_dimensional reduction_n=2}

We now investigate what SOHB system in dimension~$2$ is obtained by this dimension reduction method. With this aim, we apply Prop.~\ref{prop_dim_reduc} with~$n=3$ and~$p=2$. If we do so, we end up with the following system for~$(\rho, \Omega)$:~${\mathbb R}^2 \times [0,T] \to (0,\infty) \times {\mathbb S}^2$:
\begin{align}
&
\partial_t \rho + \nabla \cdot (c_1 \rho \Omega) = 0, \label{eq:mass_dim2} \\
&
\rho ( \partial_t + (c_2-c_4) \Omega \cdot \nabla) \Omega = - c_3 P_{\Omega^\bot} \nabla \rho + c_4 \rho \,  \textrm{curl} \, \Omega \, \, \Omega^\bot .  \label{eq:orient_dim2}
\end{align}
Letting~$\Omega_x$ and~$\Omega_y$ be the two components of~$\Omega$ i.e.~$\Omega = (\Omega_x,\Omega_y)^T$, we define~$\Omega^\bot = (-\Omega_y, \Omega_x)^T$ and~$\textrm{curl} \, \Omega$, the scalar given by 
\[ \textrm{curl} \, \Omega = \partial_{x} \Omega_y - \partial_y \Omega_x.\]
Eq.~\eqref{eq:orient_dim2} directly follows from~\eqref{eq:om1} by noticing that~$\nabla \cdot \Omega_2 = \nabla \cdot \Omega^\bot = - \textrm{curl} \, \Omega$. 

We now show that~\eqref{eq:orient_dim2} can be reduced to an equation of the form~\eqref{eq:orient_SOH} but with different coefficients. Since~$\Omega$ is a unit vector,~$(\Omega\cdot\nabla)\Omega$ is proportional to~$\Omega^\bot$. Furthermore, from~$\partial_x(\Omega_x^2+\Omega_y^2)=\partial_y(\Omega_x^2+\Omega_y^2)=0$, we deduce that 
$\Omega_x \partial_x \Omega_x = - \Omega_y \partial_x \Omega_y$ and~$\Omega_y \partial_y \Omega_y = - \Omega_x \partial_y \Omega_x$. Therefore, we have
\begin{align*}
\Omega^\bot\cdot(\Omega\cdot\nabla)\Omega &= - \Omega_y(\Omega_x\partial_x\Omega_x + \Omega_y\partial_y\Omega_x) + \Omega_x(\Omega_x\partial_x\Omega_y + \Omega_y\partial_y\Omega_y)\\
&=\Omega_{y}^2 (\partial_x\Omega_y - \partial_y\Omega_x) + \Omega_x^2 (\partial_x\Omega_y - \partial_y\Omega_x) \\
&=\textrm{curl} \,\Omega.
\end{align*}
We deduce that
\begin{equation}
(\Omega\cdot\nabla)\Omega=\textrm{curl}\, \Omega \,\Omega^\bot, 
\label{eq:curlOmperp}
\end{equation}
and~\eqref{eq:orient_dim2} simplifies into 
\[
\rho ( \partial_t + (c_2-2c_4) \Omega \cdot \nabla) \Omega = - c_3 P_{\Omega^\bot} \nabla \rho,   
\]
which is the SOH model with coefficient~$c_2$ replaced by~$c_2 - 2 c_4$.
Thus, the SOHB system in dimension 2 coincides with the SOH system but its coefficients depend on what derivation route is taken.

Interestingly, the form~\eqref{eq:mass_dim2},~\eqref{eq:orient_dim2} lends itself to other relaxation approximations of the SOH model than~\eqref{eq:mass_relax},~\eqref{eq:orient_relax}. One possibility is to directly apply Prop.~\ref{prop:relax}. If we do so, we obtain the following system: 
\begin{align}
&
\partial_t \rho + c_1 \nabla \cdot (\rho v_1) = 0, 
\label{eq:mass_relax_2_dim2} \\
&
                                  \partial_t (\rho v_1) + (c_2 - 2 c_4) \nabla \cdot (\rho v_1 \otimes v_1) + 2 (c_3-c_4) \nabla \rho \nonumber \\
  &\hspace{8cm}= - \frac{1}{\alpha} \big( (|v_1|^2 - 1) v_1 + (v_1 \cdot v_2) v_2 \big). 
\label{eq:v1_dim2} \\
&
\partial_t (\rho v_2) - 2 c_4 \nabla \cdot (\rho v_2 \otimes v_1) + c_2 \nabla \cdot (\rho v_1 \otimes v_2) 
= - \frac{1}{\alpha} \big( (v_1 \cdot v_2) v_1 + (|v_2|^2 - 1) v_2 \big).
\label{eq:v2_dim2}
\end{align}

However, we can also directly check that the following is also a conservative relaxation approximation of System~\eqref{eq:mass_dim2} ,~\eqref{eq:orient_dim2}: 
\begin{align}
&
\partial_t \rho + \nabla \cdot (\rho v) = 0, \label{eq:mass_relax_SOHB} \\
&
\partial_t (\rho v) + \frac{1}{c_1} \nabla \cdot \big( \rho ( (c_2-c_4) v \otimes v + c_4 v^\bot \otimes v^\bot) \big) + c_1 (c_3 - c_4) \nabla \rho  \nonumber \\ 
&\hspace{10cm}
= \frac{1}{\alpha} \rho v (c_1^2 - |v|^2),  \label{eq:orient_relax_SOHB}
\end{align} 
Indeed, applying the projection~$P_{v^\bot} = \frac{1}{|v|^2} v^\bot \otimes v^\bot$ to~\eqref{eq:orient_relax_SOHB} leads to 
\begin{align} 
&
\rho P_{v^\bot} \partial_t v + \frac{(c_2-c_4)}{c_1} \rho P_{v^\bot} (v\cdot \nabla) v + \frac{c_4}{c_1} \big( (v^\bot \cdot \nabla \rho) v^\bot + \rho (\nabla \cdot v^\bot) v^\bot \big) \nonumber \\
&\hspace{10cm}
+ c_1 (c_3 - c_4) P_{v^\bot} \nabla \rho = 0 .
\label{eq:orient_relax_limit}
\end{align}
Now, from~\eqref{eq:orient_relax_SOHB}, we get that~$v \to c_1 \Omega$ with~$|\Omega|=1$ as~$\alpha \to 0$. So, taking the limit~$\alpha \to 0$ in~\eqref{eq:orient_relax_limit}, we get 
\[
\rho \partial_t \Omega + (c_2-c_4) \rho (\Omega \cdot \nabla) \Omega + c_4 \big( (\Omega^\bot \cdot \nabla \rho) \Omega^\bot + \rho (\nabla \cdot \Omega^\bot) \Omega^\bot \big) + (c_3 - c_4) P_{\Omega^\bot} \nabla \rho = 0, 
\]
which, owing to the fact that~$(\Omega^\bot \cdot \nabla \rho) \Omega^\bot = P_{\Omega^\bot} \nabla \rho$ and~$\nabla \cdot \Omega^\bot = - \textrm{curl} \, \Omega$, leads to~\eqref{eq:orient_dim2}.

The relaxation systems~\eqref{eq:mass_relax_2_dim2}-\eqref{eq:v2_dim2} on the one hand and~\eqref{eq:mass_relax_SOHB},~\eqref{eq:orient_relax_SOHB} on the other hand do not coincide. Indeed, in the latter,~$v_2 = v^\bot$ is constrained to be orthogonal to~$v_1 = v$ while in the former, there is no a priori relation between~$v_1$ and~$v_2$, both being subject to evolution equations in their own right. This shows that there may exist several distinct conservative relaxation approximations of the SOHB model. These different relaxation approximations may not be equivalent in terms of numerical efficiency. Further studies are needed to select the most appropriate ones. 

Because of~\eqref{eq:curlOmperp}, we note that if we change one sign in~\eqref{eq:orient_dim2}, namely writing 
\begin{equation}
\rho ( \partial_t + (c_2+c_4) \Omega \cdot \nabla) \Omega = - c_3 P_{\Omega^\bot} \nabla \rho + c_4 \rho \,  \textrm{curl} \, \Omega \, \, \Omega^\bot ,  
\label{eq:orient_SOH_bis} 
\end{equation}
we obtain an equation equivalent to the original SOH orientation equation~\eqref{eq:orient_SOH}, whatever the value of~$c_4$ is. The previous discussion shows that the following is a relaxation approximation of~\eqref{eq:orient_SOH_bis}: 
\[\partial_t (\rho v) + \frac{1}{c_1} \nabla \cdot \big( \rho ( (c_2+c_4) v \otimes v + c_4 v^\bot \otimes v^\bot) \big) + c_1 (c_3 - c_4) \nabla \rho 
= \frac{1}{\alpha} \rho v (c_1^2 - |v|^2),\]
which can also be rewritten, using that~$v\otimes v + v^\bot \otimes v^\bot = |v|^2 \textrm{I}~$, as
\begin{equation}
\partial_t (\rho v) + \frac{c_2}{c_1} \nabla \cdot \big( \rho v \otimes v \big) + c_1 c_3 \nabla \rho - \frac{c_4}{c_1} \nabla \big( \left(c_1^2-|v|^2\right) \rho \big)  
= \frac{1}{\alpha} \rho v (c_1^2 - |v|^2). \label{eq:relax_SOH_c4}
\end{equation}

Two equations~\eqref{eq:relax_SOH_c4} for different values of~$c_4$ are not equivalent (i.e. do not generate the same solution even if started with the same initial condition). Consequently, we obtain a one-parameter family of non-equivalent relaxation approximations of the SOH system 
\eqref{eq:mass_SOH}-\eqref{eq:orient_SOH} of which System~\eqref{eq:mass_relax},~\eqref{eq:orient_relax} is the special case~$c_4 = 0$. This again shows that there might exist many (and indeed, uncountably many) different relaxation approximations of the SOH model.

\setcounter{equation}{0}
\section{Non-conservativity: second approach by viscosity}
\label{sec:viscosity}

As exposed in Section~\ref{subsec_relax_rationale}, discontinuous solutions (aka shock waves) of
non-conservative hyperbolic models are not uniquely defined. Microscopic dissipation mechanisms which are needed to single out the ``right'' solution are lost in the macroscopic limit. One way out of this problem is to modify the scaling assumptions so as to maintain these dissipation mechanisms finite in the macroscopic limit. The goal of the present section is to develop this approach and derive a viscous SOHB model. 

In~\cite{degond2021body}, the SOHB model was derived from a kinetic model describing a system of swarming rigid bodies. We return to this framework but modify the background scaling in order to retrieve viscous correction to the SOHB model derived in~\cite{degond2021body}.

\subsection{The kinetic framework}
\label{subsec_kin_frame}

The object under study is the the probability distribution function~$f(x,A,t) \, dx \, dA$ of particles in the small volume~$dx \, dA$ about position~$x \in {\mathbb R}^n$ and rotation~$A \in \textrm{SO}_n {\mathbb R}$ at time~$t \in [0,T]$, where~$A$ represents the rotation mapping the reference frame~$(e_1, \ldots, e_n)$ to the frame attached to the corresponding rigid body. The quantity~$f$ is subject to 
\begin{equation}
\partial_t f + (c_0 A e_1 \cdot \nabla_x) f = \nu \,  \big( \rho_f M_{\tilde \rot_f} - f \big) . 
\label{eq:BGK_unsc}
\end{equation} 
where~$c_0$ is the particle speed (which is a fixed constant),~$A e_1$ is the direction of self-propulsion of the particles with attached frame represented by~$A$, and~$\rho_f$:~${\mathbb R}^n \times [0,T] \to [0,\infty)$ is the local density 
\begin{equation} 
\rho_f (x,t) = \int_{\mathrm{SO}_n {\mathbb R}} \, f(x, A,t) \, dA, 
\label{eq:def_rho}
\end{equation}
with~$dA$ being the Haar measure on~$\mathrm{SO}_n {\mathbb R}$. The function ~$\tilde \rot_f$:~${\mathbb R}^n \times [0,T] \to \mathrm{SO}_n {\mathbb R}$ is defined by  
\begin{equation*}
\tilde \rot_f(x,t) = {\mathcal P} (\tilde J_f(x,t)), 
\end{equation*}
where, for an~$n \times n$ matrix~$J$ such that~$\det J \not = 0$,~${\mathcal P} (J)$ is the unique element of~$\mathrm{SO}_n {\mathbb R}$ such that 
\begin{equation*}
{\mathcal P} (J) := \textrm{arg  max}_{A \in \mathrm{SO}_n {\mathbb R}} \, A \cdot J , 
\end{equation*}
and~$\tilde J_f(x,t)$ is defined by 
\begin{equation}
\tilde J_f(x,t) =  \int_{{\mathbb R}^n \times \mathrm{SO}_n {\mathbb R}} K(x-y) \, f(y, A,t) \, A \, dy \, dA. 
\label{eq:def_tilJ}
\end{equation}
In~\eqref{eq:def_tilJ},~$K$ is a given sensing function:~${\mathbb R}^n \to [0,\infty)$ supposed to be radial, i.e.  
\begin{equation} 
K(x) = \frac{1}{R^n} \tilde K(\frac{|x|}{R}), 
\label{eq:sensing}
\end{equation}
for a given~$\tilde K$:~$[0,\infty) \to [0,\infty)$ and a spatial scale parameter~$R>0$ referred to as the sensing radius. We will assume that the condition~$\det \tilde J_f(x,t) > 0$ is satisfied everywhere. If~$\det J > 0$,~${\mathcal P} (J)$ is nothing but the orthogonal matrix in the polar decomposition of~$J$. Somehow,~$\tilde \rot_f$ represents the average body attitude of the particles in the neighborhood of~$x$. The parameter~$\nu > 0$ is a fixed relaxation rate. For a given~$\rot \in \textrm{SO}_n {\mathbb R}$,~$M_\rot$:~$\textrm{SO}_n \to {\mathbb R}$, denotes the von Mises distribution centered at~$\rot$ given by  
\begin{equation*}
M_{\rot} (A) = \frac{1}{Z} \exp (\kappa \, \rot \cdot A), \quad Z = \int_{\mathrm{SO}_n {\mathbb R}} \exp \big( \kappa \, \mbox{Tr} A \big)  \, dA,
\end{equation*}
where~$\kappa >0$ is a given constant playing the role of an inverse temperature and called the concentration parameter. The rotation~$\rot$ plays the role of the mean of the von Mises distribution and is referred to as the orientation parameter. Eq.~\eqref{eq:BGK_unsc} describes a system of rigid bodies, which are self-propelled at speed~$c_0$ in the direction~$A e_1$ and which mutually interact at rate~$\nu$. Through their interaction, the distribution function relaxes to a von-Mises distribution with orientation parameter equal to the average body attitude~$\tilde \rot_f$ in the neighborhood of~$x$. In other words, during a time~$dt$ a proportion~$\nu dt$ of the particles located in the neighborhood of~$x$ tend to adopt a body attitude close to~$\tilde \rot_f$ up to a small noise measured by~$\kappa^{-1}$.  

We now scale the kinetic equation~\eqref{eq:BGK_unsc}. Let~$t_0$ and~$x_0 = c_0 t_0$ be a choice of time and space scales. We introduce dimensionless variables~$t' = t/t_0$ and~$x' = x/x_0$ and scale the functions~$f$,~$\rho_f$ and~$\tilde J_f$ as follows: 
\[ f'(x',A,t') = x_0^n f(x,A,t), \quad \rho_{f'}'(x',t') = x_0^n \rho_f(x,t),  \quad \tilde J_{f'}'(x',t') = x_0^n \tilde J_f(x,t).\] 
Through this scaling,~\eqref{eq:BGK_unsc} is transformed into 
\begin{equation}
\partial_t f + (A e_1 \cdot \nabla_x) f = \bar \nu \,  \big( \rho_f M_{\tilde \rot_f} - f \big) , 
\label{eq:BGK_sc}
\end{equation} 
where the primes have been dropped for simplicity and where~$\bar \nu = \nu t_0$ is a dimensionless relaxation frequency. Additionally, thanks to~\eqref{eq:sensing}, we have 
\begin{equation} 
\tilde J_f(x) =   \int_{{\mathbb R}^n \times \mathrm{SO}_n {\mathbb R}} \frac{1}{\bar R^n} \tilde K \big(\frac{|x-y|}{\bar R} \Big) \, f(y, A) \, A \, dy \, dA, 
\label{eq:tilJ}
\end{equation}
and~$\bar R = R/x_0$ is the dimensionless sensing radius. We now make the following fundamental scaling assumptions:
\begin{equation} 
\frac{1}{\bar \nu} = \varepsilon \ll 1, \qquad \bar R = \sqrt{\varepsilon} ,
\label{eq:scaling_ass}
\end{equation}
where~$\varepsilon >0$ is the dimensionless rate at which~$f$ relaxes towards the local equilibrium~$\rho_f M_{\rot_f}$. Eq.~\eqref{eq:scaling_ass} assumes that the sensing radius is small compared with the macroscopic length scale (as~$\sqrt{\varepsilon} \ll 1$) and large compared with the microscopic length scale (as~$\sqrt{\varepsilon} \gg \varepsilon$). This scaling was introduced in~\cite{degond20133hydrodynamic} and referred to as the weakly non-local asymptotics. A similar scaling was introduced in~\cite{wang2015small, zhang2006molecular} in the Doi-Onsager equation to recover elastic stresses in the small Deborah number limit. This is where the present derivation departs from previous works~\cite{degond2021body, degondfrouvellemerino17, degondfrouvellemerinotrescases18, degond2018alignment} which assumed~$\bar R = \varepsilon$. The weakly non-local scaling is responsible for the appearance of diffusive terms in the macroscopic equations which are not present if the classical scaling is used.   

Taylor expanding~\eqref{eq:tilJ} (with~$R =\sqrt{\varepsilon}$) with respect to~$\varepsilon$ and using rotational symmetry, we get 
\[ \tilde J_f  = J_f + \varepsilon \alpha \Delta J_f + {\mathcal O}(\varepsilon^2), \]
with
\begin{equation}
 J_f(x) =  \int_{\mathrm{SO}_n {\mathbb R}}  \, f(x, A) \, A \, dA \quad \textrm{ and } \quad \alpha = \frac{1}{2n} \int_{{\mathbb R}^n} \tilde K(|\xi|) \, |\xi|^2 \, d \xi. 
\label{eq:def_Jf}
\end{equation}
Now, denoting by~${\mathcal M}^+$ the space of~$n \times n$ matrices of positive determinant, the map 
${\mathcal P}$:~${\mathcal M}^+ \to \textrm{SO}_n{\mathbb R}$,~$J \mapsto {\mathcal P} (J)$ maps~$J$ to the orthogonal factor of the polar decomposition of~$J$. Thus, this map is differentiable on~${\mathcal M}^+$ and we denote by~$d_J {\mathcal P} (J)$ its differential at~$J$ which is a linear map from~${\mathcal M}$ to~$T_{{\mathcal P} (J)}$, the tangent space to~$\textrm{SO}_n{\mathbb R}$ at~${\mathcal P} (J)$. An explicit expression of~$d_J {\mathcal P} (J)$ is not needed at this stage. Thus, we have 
\[ \tilde \rot_f = \rot_f + \varepsilon \alpha d_J{\mathcal P}(J_f) (\Delta J_f) + {\mathcal O}(\varepsilon^2),\]
with~$ \rot_f = {\mathcal P}(J_f)$. Now, for a given~$A \in \textrm{SO}_n {\mathbb R}$, the map~$\textrm{SO}_n {\mathbb R} \to {\mathbb R}$,~$\rot \mapsto M_{\rot}(A)$ has differential~$d_{\rot} M_{\rot}(A)$. This is a linear map from~$T_{\rot}$ to~${\mathbb R}$ given by 
\[ d_{\rot} M_{\rot}(A) (Q) = \kappa M_{\rot}(A) P_{T_\rot} A \cdot Q, \qquad \forall Q \in T_\rot.\]
We refer to the proof of Prop.~\ref{prop:relax} for the expressions of~$T_\rot$ and~$P_{T_\rot}$. By the chain rule, it follows that 
\[ M_{\tilde \rot_f} = M_{\rot_f} + \varepsilon \alpha \kappa \, M_{\rot_f}(A) \, P_{T_{\rot_f}} A \cdot \big( d_J{\mathcal P}(J_f) (\Delta J_f) \big) + {\mathcal O}(\varepsilon^2).\]
Inserting this expansion in~\eqref{eq:BGK_sc}, we get 
\begin{align}
&
\partial_t f^\varepsilon + (A e_1 \cdot \nabla_x) f^\varepsilon - \alpha \kappa \, \rho_{f^\varepsilon} M_{\rot_{f^\varepsilon}}(A) \, P_{T_{\rot_{f^\varepsilon}}} A \cdot \big( d_J{\mathcal P}(J_{f^\varepsilon}) (\Delta J_{f^\varepsilon}) \big) \nonumber \\
&\hspace{10cm}
= \frac{1}{\varepsilon} \,  \big( \rho_{f^\varepsilon} M_{\rot_{f^\varepsilon}} - f^\varepsilon \big) , 
\label{eq:BGK_eps}
\end{align}
where~$\rho_f$ and~$J_f$ are defined by~\eqref{eq:def_rho} and~\eqref{eq:def_Jf} respectively, and~$\rot_f = {\mathcal P}(J_f)$. In~\eqref{eq:BGK_eps}, we have highlighted the dependence of~$f$ on~$\varepsilon$ and we have neglected the~${\mathcal O}(\varepsilon)$ terms. Indeed, it may be checked that these higher order terms have no contribution to the final result. In the next section, we aim to find the formal limit~$\varepsilon \to 0$ of this problem.

\subsection{Derivation of the viscous SOHB model}
\label{subsec_visc_sohb}

The limit~$\varepsilon \to 0$ of~\eqref{eq:BGK_eps} is formalized in the following

\begin{theorem}
Assume~$n \geq 3$. Let~$f^\varepsilon$:~${\mathbb R}^n \times \textrm{SO}_n{\mathbb R} \times [0,T] \to {\mathbb R}$ be a solution to~\eqref{eq:BGK_eps}. Assume that it is smooth and converges to a smooth~$f$ defined on the same set as~$\varepsilon \to 0$. Then,~$f= \rho(x,t) M_{\rot(x,t)}(A)$, where~$(\rho,\rot)$: ~${\mathbb R}^n \times [0,T] \to [0,\infty) \times \textrm{SO}_n{\mathbb R}$ is a solution to 
\begin{align}
&
\partial_t \rho + \nabla \cdot (c_1 \rho \Omega_1) = 0,  \label{eq:mass_visc} \\
&
\rho \big( \partial_t + c_2 \, \Omega_1 \cdot \nabla \big) \rot - {\mathbb W} \rot = \alpha P_{T_\rot} \Delta (\rho \rot),  
\label{eq:orient_visc} 
\end{align}
where~$\Omega_1$ and~${\mathbb W}$ are respectively given by~\eqref{eq:def_Omega1},~\eqref{eq:expressF} and~$P_\rot$ by~\eqref{eq:def_Ptheta}. The constants~$c_1$, \ldots~$c_4$ involved in the model are functions of~$\kappa$ and~$n$. Their expressions are given by formulas (3.11)-(3.14) of~\cite{degond2021body}. 
\label{thm:deriv_visc}
\end{theorem}

\medskip
\noindent
We note that the continuity equation~\eqref{eq:mass_visc} is identical with that appearing in the non-viscous case~\eqref{eq:mass_2}. This is a classical feature of viscous fluid models. We also remark that in dimension~$n=2$, we would find the same result, except for the coefficient~$c_4$ in the expression of~${\mathbb W}$ being equal to~$0$ (see beginning of Section~\ref{subsec_relax_dim2}). 

\medskip
\noindent
\textbf{Proof.} In the absence of the third term at the left-hand side of~\eqref{eq:BGK_eps} and in the case~$n \geq 3$, this theorem has been proved in~\cite{degond2021body}. So, we will only focus on the treatment of this additional term. Denote by 
\begin{align*}
&
T_1 f =: \partial_t f + (A e_1 \cdot \nabla_x) f, \\
&
T_2 f =: - \alpha \kappa \, \rho_f M_{\rot_f}(A) \, P_{T_{\rot_f}} A \cdot \big( d_J{\mathcal P}(J_f) (\Delta J_f) \big) , 
\end{align*}
It has been shown (see formulas (4.2) and (4.7) of~\cite{degond2021body}) that~\eqref{eq:mass_2} and~\eqref{eq:orient_visc} are respectively equivalent to 
\begin{align}
&
\int_{\textrm{SO}_n{\mathbb R}} (T_1 + T_2) (\rho M_\rot) \, dA = 0, \label{eq:intdA}\\
&
\int_{\textrm{SO}_n{\mathbb R}} (T_1 + T_2) (\rho M_\rot) \, P_{T_\rot} A \, dA = 0. \label{eq:intPTAdA}
\end{align}
Moreover, the contribution of~$T_1 (\rho M_\rot)$ to~\eqref{eq:intdA} leads to the left hand side of~\eqref{eq:mass_2} while its contribution to~\eqref{eq:intPTAdA} leads to the left hand side of~\eqref{eq:orient_3} multiplied by~$\kappa C_2$, where~$C_2$ is given by formula (4.14) of~\cite{degond2021body}. 

We first compute~$T_2 (\rho M_\rot)$. We have 
\[ J_{\rho M_\rot} = \rho \int_{\textrm{SO}_n{\mathbb R}} M_\rot(A) \, A \, dA = \rho c_1 \rot,\]
by Lemma 4.1 of~\cite{degond2021body}. So,~$J_{\rho M_\rot}$ is of the form~$\lambda \rot$ for~$\lambda >0$ and we need to compute~$d_J{\mathcal P}(\lambda \rot)$. We write~$J = \lambda \rot + \delta J_1$ with~$\delta \ll 1$ and evaluate~${\mathcal P}(\lambda \rot + \delta J_1)$ at first order in~$\delta$. Since~$\det (\lambda \rot) = \lambda^n >0$,~$\lambda \rot \in {\mathcal M}^+$ and for~$\delta$ small enough, then~$\lambda \rot + \delta J_1 \in {\mathcal M}^+$. Then,~${\mathcal P}(\lambda \rot + \delta J_1)$ is nothing but the orthogonal matrix in the polar decomposition of~$\lambda \rot + \delta J_1$. We recall that we have the following expression:  
\[{\mathcal P} (J) = (J J^T)^{-1/2} J, \quad \forall J \in {\mathcal M}^+, \]
where~$(J J^T)^{-1/2}$ is the inverse of the matrix~$(J J^T)^{1/2}$ and the definition of the square root of a symmetric positive definite matrix has been recalled in Section~\ref{subsec_relax_dimn}. So, 
\begin{align} 
{\mathcal P}(\lambda \rot + \delta J_1) &= \big[ (\lambda \rot + \delta J_1)(\lambda \rot^T + \delta J_1^T)\big]^{-1/2} (\lambda \rot + \delta J_1) \nonumber \\
&= \big[ (\rot + \tilde \delta J_1)(\rot^T + \tilde \delta J_1^T)\big]^{-1/2} (\rot + \tilde \delta J_1) \nonumber \\
&= \big[ \textrm{I} + \tilde \delta (J_1 \rot^T + \rot J_1^T) + {\mathcal O}(\tilde \delta^2) \big]^{-1/2} (\rot + \tilde \delta J_1), \label{eq:Plrot+delJ1}
\end{align}
with~$\tilde \delta = \delta / \lambda$. In a neighbourhood of~$\textrm{I}$, we have the normally convergent series: 
\[ (\textrm{I} + A)^{-1/2} = \sum_{k=0}^\infty \left( \begin{array}{c} - 1/2 \\ k \end{array} \right) A^k = \textrm{I} - \frac{1}{2} A + {\mathcal O}(\|A\|^2),\]
where~$\| \cdot \|$ is any matrix norm, such as the Frobenius norm and we have used the generalized choose symbol.  We can insert this expansion in~\eqref{eq:Plrot+delJ1} and finally get 
\[ {\mathcal P}(\lambda \rot + \delta J_1) = \rot + \frac{\delta}{\lambda} P_{T_\rot} J_1 +  {\mathcal O}(\tilde \delta^2),\]
from which we find that 
\[d_J{\mathcal P}(\lambda \rot)(J_1) = \frac{1}{\lambda} P_{T_\rot} J_1.\]
We easily deduce that 
\begin{equation}
T_2 (\rho M_\rot) = - \alpha \kappa M_\rot(A) \, P_{T_\rot} A \cdot P_{T_\rot} \big(\Delta (\rho \rot) \big). 
\label{eq:T2rhoM}
\end{equation}

We first remark that 
\[ \int_{\textrm{SO}_n{\mathbb R}} M_\rot(A) \, P_{T_\rot} A \, dA = P_{T_\rot} \Big( \int_{\textrm{SO}_n{\mathbb R}} M_\rot(A) \, A \, dA \Big)  =  c_1 P_{T_\rot} \rot = 0,\]
from which we infer that 
\begin{equation} 
\int_{\textrm{SO}_n{\mathbb R}} T_2 (\rho M_\rot) \, dA = 0 , 
\label{eq:intT2}
\end{equation}
showing that the continuity equation~\eqref{eq:mass_2} is unmodified by the presence of the extra term~$T_2$ (note that~\eqref{eq:intT2} would still be true with~$\rho M_\rot$ replaced by any distribution function~$f$). 

We now consider the term 
\begin{equation} 
{\mathcal T}_2 = \int_{\textrm{SO}_n{\mathbb R}} T_2 (\rho M_\rot) \, P_{T_\rot} A \, dA . 
\label{eq:intT2rhoM}
\end{equation}
Inserting~\eqref{eq:T2rhoM} into~\eqref{eq:intT2rhoM} and using the notations~$M = \Delta (\rho \rot)$ and~$P = (M \rot^T - \rot M^T)/2$ leads to 
\begin{align*} {\mathcal T}_2 &= - \alpha \kappa \int_{\textrm{SO}_n{\mathbb R}} \Big( \frac{A \rot^T - \rot A^T}{2} \rot \cdot \frac{M \rot^T - \rot M^T}{2} \rot \Big) \frac{A \rot^T - \rot A^T}{2} \rot \, M_\rot(A) \, dA  \\
&= - \alpha \kappa \int_{\textrm{SO}_n{\mathbb R}} \Big( \frac{A - A^T}{2} \cdot P \Big) \frac{A - A^T}{2}  \, M_{\textrm{I}}(A) \, dA \, \, \rot \\
&= - \alpha \kappa \int_{\textrm{SO}_n{\mathbb R}}(A \cdot P) \frac{A - A^T}{2}  \, M_{\textrm{I}}(A) \, dA \, \, \rot, 
\end{align*}
where we changed variable from~$A$ to~$A \rot^T$ and exploited the translation invariance of the Haar measure, and where we used that~$P$ is antisymmetric. In~\cite{degond2021body}, it was proved that the integral appearing in the last line equals~$C_2 P$ where~$C_2$ is given by formula~(4.14) of~\cite{degond2021body}. Thus, 
\begin{equation} 
{\mathcal T}_2 = - C_2 \kappa \alpha P \rot = - C_2 \kappa \alpha P_{T_\rot} \Delta (\rho \rot). 
\label{eq:calT2}
\end{equation}
According to the discussion at the beginning of the proof, we eventually find that the left-hand side of~\eqref{eq:orient_3} must be complemented with~\eqref{eq:calT2} divided by~$C_2 \kappa$. This results in~\eqref{eq:orient_visc} and ends the proof.  \endproof

\setcounter{equation}{0}
\section{Conclusion and perspectives}
\label{sec:conclu}

In this paper, we have studied the SOHB model which provides a hydrodynamic description of a system of swarming rigid bodies interacting through alignment. We have shown that the model is hyperbolic and that it can be approximated by a conservative system through a relaxation limit. We have also derived viscous corrections to the SOHB model from the kinetic model thanks to a small variant of the scaling assumptions compared to earlier work~\cite{degond2021body}. The goal of this study is to provide the bases for numerical approximations of the SOHB model. Indeed, its nonconservative character requires special treatment. One possibility is to use a splitting method previously developed for the SOH model which centrally uses a conservative relaxation approximation. Another possibility is to take advantage of the viscous terms to prevent the appearance of discontinuous solutions. These two approaches will be explored in future work.



\newpage
\appendix

\setcounter{equation}{0}
\section{Invariance by coordinate change: direct proof of Prop.~\ref{prop:invariance_coordinate_change}}
\label{sec_invar_coor_change_proof}

Let~$(f'_1, \dots, f'_n)$ be another direct orthonormal frame and~$x'=(x'_1, \ldots, x'_n)^T$ the coordinates of the same point~${\mathbf x}$ in this new frame, i.e.~${\mathbf x} = \sum_{i=1}^n x'_i \, f'_i$. Let~$S$ be the transition matrix from~$(f_1, \ldots, f_n)$ to~$(f'_1, \ldots, f'_n)$. Then,~$x = S x'$. If~$\rot$ is the matrix of the map that sends~$(e_1, \ldots, e_n)$ to~$(\Omega_1, \ldots, \Omega_n)$ expressed in the basis~$(f_1, \ldots, f_n)$, then the matrix of the same map expressed in the basis~$(f'_1, \dots, f'_n)$ is~$\rot' = S^T \rot S$. Thus, the change of coordinate frame~$(f_1, \ldots, f_n)$ to~$(f'_1, \dots, f'_n)$ leads to the following change of variables and unknowns: 
\begin{equation*}
\rho'(x',t) = \rho(Sx',t), \quad \rot'(x',t) = S^T \rot(Sx',t) S.
\end{equation*}

We can directly prove this by remarking that if~$P$ and~$Q$ are the transition matrices from~$(f_1, \ldots, f_n)$ to~$(e_1, \ldots, e_n)$ and~$(\Omega_1, \ldots, \Omega_n)$ respectively, then,~$\rot = Q P^T$. Indeed, by definition,~$\rot$ is the transition matrix from~$(f_1, \ldots, f_n)$ to~$(\rot f_1, \ldots, \rot f_n)$. But, by linearity,~$P$ is also the transition matrix from~$(\rot f_1, \ldots, \rot f_n)$ to~$(\Omega_1, \ldots, \Omega_n)$. Thus, the transition matrix from~$(f_1, \ldots, f_n)$ to~$(\rot f_1, \ldots, \rot f_n)$ is equal to~$QP^T$, showing the claim. Now, if~$P'$ and~$Q'$ are the transition matrices from~$(f'_1, \ldots, f'_n)$ to~$(e_1, \ldots, e_n)$ and~$(\Omega_1, \ldots, \Omega_n)$ respectively, we have~$\rot' = Q' P^{'T}$. But~$P=SP'$ and~$Q = SQ'$ so that~$\rot' = S^T Q P^T S = S^T \rot S$ as stated. 

Now, we need to point out a slight abuse of notation which was used so far. When we said~$\Omega_k = \rot e_k$, we should have rather said~$\omega_k = \rot \epsilon_k$ where~$\omega_k$ and~$\epsilon_k$ are the~$n$-tuple of coordinates of~$e_k$ and~$\Omega_k$ in the coordinate system~$(f_1, \ldots, f_n)$. Note that the present~$\omega_k$ has no relation with that appearing in Section~\ref{sec_hyperbolic}. These same vectors~$e_k$ and~$\Omega_k$ have~$n$-tuple of coordinates~$\epsilon_k'$ and~$\omega_k'$ in the coordinate system~$(f'_1, \ldots, f'_n)$ with 
\begin{equation}
\omega_k = S \omega_k' \quad \textrm{ and } \quad \epsilon_k = S \epsilon'_k.
\label{eq:om_eps}
\end{equation}
Thus, we get
\[ \omega_k'(x',t) = S^T \omega_k(Sx',t) = S^T \rot(Sx',t) \epsilon_k = S^T \rot(Sx',t) S \epsilon'_k = \rot'(x',t) \epsilon_k',\]
which, in particular, shows that the relation~$\Omega_1 = \rot e_1$ (with our abuse of notation) is invariant by coordinate change. 

Now, the vector~$\rho \Omega_1$ has coordinates~$\rho(x,t) \omega_1(x,t)$ and~$\rho'(x',t) \omega'_1(x',t)$ in the two coordinate frames with~$ (\rho' \omega'_1)(x',t) = (\rho S^T \omega_1)(Sx',t)$. Thus, 
\begin{align}
 \nabla_{x'} \cdot (\rho' \omega_1') (x',t) &= \sum_{i=1}^n \partial_{x'_i} \big[ (\rho (S^T \omega_1)_i) (Sx',t)\big] = \sum_{i,j,k = 1}^n S_{ji} \partial_{x_k}(\rho (\omega_1)_j) (Sx',t) S_{ki} \nonumber \\
&= \sum_{j=1}^n \partial_{x_j} (\rho (\omega_1)_j) (Sx',t) = \nabla_x \cdot (\rho \omega_1) (Sx',t). 
\label{eq:div_coord_change}
\end{align}
The third equality is due to the fact that~$\sum_{i=1}^n S_{ji} S_{ki} = \delta_{jk}$ where~$\delta_{jk}$ is the Kronecker delta and this identity comes from~$S$ being an orthogonal matrix. From this, it follows that~\eqref{eq:mass_2} is invariant by coordinate change. 

Now, for any scalar function~$\varphi(x)$, letting~$\varphi'(x') = \varphi(Sx')$, we have 
\begin{equation}
\nabla_{x'} \varphi'(x') = S^T \nabla_x \varphi(Sx'). 
\label{eq:grad_coord_change}
\end{equation}
It follows that 
\[ \omega'_1 (x') \cdot \nabla_{x'} \varphi'(x') = (S^T \omega_1(Sx')) \cdot (S^T \nabla_x \varphi(Sx')) = (\omega_1 \cdot \nabla_x \varphi) (Sx'),\]
again using the fact that~$S$ is an orthogonal matrix. Hence, the operator~$\Omega_1 \cdot \nabla$ (still with our abuse of notation) and consequently, the convective derivative~$\partial_t + c_2 \Omega_1 \cdot \nabla$ are invariant by coordinate change. 

Now, using~\eqref{eq:Thet_nab_Thet} and~\eqref{eq:om_eps} as well as a similar computation as in~\eqref{eq:div_coord_change}, we get that 
\[ (\rho' \rot' \nabla_{x'} \cdot \rot')(x',t) = S^T (\rho \rot \nabla_x \cdot \rot)(Sx',t).\]
This combined with~\eqref{eq:grad_coord_change} shows that the coordinates~${\mathcal F}(x,t)$ and~${\mathcal F}'(x',t)$ of~$F$ in the two coordinate systems are related by the same relation as~\eqref{eq:om_eps}. We deduce that 
\[({\mathcal F}' \wedge \omega_1')(x',t) = S^T ({\mathcal F} \wedge \omega_1)(Sx',t) S.\]
Now, similar computations as those done in~\eqref{eq:div_coord_change} and left to the reader show that 
\[ (\nabla_{x'} \wedge \omega'_1) (x',t) = S^T  (\nabla_x \wedge \omega_1) (Sx',t) S,\]
so that~${\mathbb W}$ follows the same transformation. So, multiplying~\eqref{eq:orient_3} on the left by~$S^T$ and on the right by~$S$ and applying the formulas above shows that~\eqref{eq:orient_3} is invariant by coordinate change, which finishes the proof. \endproof

\setcounter{equation}{0}
\section{Dimension reduction: proof of Prop.~\ref{prop_dim_reduc}}
\label{sec_dimension_reduc_proof}

As noticed in Subsection~\ref{subsubsec_dimensional reduction}, the result of~Proposition~\ref{prop_dim_reduc} is a local in time property, so we only need to prove it on a small interval of time. We write~$\omega_j^0=\Pi(\Omega_j^0)$ for~$1\leq j\leq p$. By the condition~$(ii)$ of Hypothesis~\ref{hyp:dim_reduc}, these are orthogonal unit vectors which can be expressed in the basis~$(\bar e_1',\ldots,\bar e_p')$, and we have~$i(\omega_j^0)=\Omega_j^0$. By assumption, we denote by~$\bar \rho,\theta$ a solution of the SOHB system~\eqref{eq:mass_2}-\eqref{eq:expressF} in dimension~$p$ (in the frame~$(\bar e_1',\ldots,\bar e_p')$ for the rotation variable and~$(f_1,\ldots,f_p)$ for the space variable), on a small interval~$[0,t]$, with initial condition~$\bar \rho^0$ and~$\theta^0$ (given by~$\theta^0\bar e_j'=\omega_j^0$). Writing the system for the variables~$\omega_j=\theta\bar e_j'$, we get, as in~\eqref{eq:mass_rhom}-\eqref{eq:omj} :
\begin{align}
&\partial_t \bar \rho + \bar \nabla \cdot (c_1 \bar \rho \omega_1) = 0,\label{eq:mass_rhom_reduced} \\ 
& \bar \rho ( \partial_t + (c_2-c_4) \omega_1 \cdot \bar \nabla) \omega_1 = - c_3 \bar P_{\omega_1^\bot} \bar \nabla \bar \rho - c_4 \bar \rho \sum_{k=2}^p (\bar \nabla \cdot \omega_k) \omega_k , \label{eq:om1_reduced} \\
& \bar \rho ( \partial_t + c_2 \omega_1 \cdot \bar \nabla) \omega_j = \Big( c_3 \omega_j \cdot \bar \nabla \rho + c_4 \bar \rho \bar \nabla \cdot \omega_j \Big) \omega_1 \nonumber \\
&\hspace{4cm}+ c_4 \bar \rho \Big( (\bar \nabla \omega_j) \omega_1 + (\omega_j \cdot \bar \nabla) \omega_1 \Big), \quad \forall j \in \{2, \ldots, p \}, \label{eq:omj_reduced}
\end{align}
where we have introduced the notation~$\bar \nabla = (\partial_{x_1}, \ldots, \partial_{x_p})$, and~$\bar P_{\omega_1^\bot}$ is the orthogonal projection onto~$\{\omega_1\}^\bot$ in~$\textrm{Span} \{ f_1, \ldots, f_p \}$.

We now set~$\Omega_j(x,t)=i(\omega_j)(\Pi x,t)$ for~$1\leq j\leq p$,~$\Omega_j(x,t)=\Omega_j^0$, and~$\rho(x,t)=\bar \rho(\Pi x,t)$ (on~$\mathbb{R}^n\times[0,t]$).
We easily realize that~$\nabla \cdot (\rho \Omega_1) = \bar \nabla \cdot (\bar \rho \omega_1)$. Hence, from~\eqref{eq:mass_rhom_reduced} we obtain~\eqref{eq:mass_rhom}.  

Next, we apply the injection~$i$ to Eqs.~\eqref{eq:om1_reduced},~\eqref{eq:omj_reduced}. We have~$ i(\omega_1 \cdot \bar \nabla) \omega_j=(\Omega_1 \cdot \nabla) \Omega_j$ for~$1 \leq j \leq p$. Likewise,~$i(\bar \nabla \bar \rho)=\nabla \rho$ and~$ \omega_1 \cdot \bar \nabla \bar \rho=\Omega_1 \cdot \nabla \rho$ as both~$\Omega_1$ and~$\nabla \rho$ lie in~$\textrm{Span} \{ f_1, \ldots, f_p \}$. Thus,~$i(\bar P_{\omega_1^\bot} \bar \nabla \bar \rho) = P_{\Omega_1^\bot} \nabla \rho$. We also have~$\nabla \cdot \Omega_k = \bar \nabla \cdot \omega_k$ for all~$1\leq k \leq p$, while~$\nabla \cdot \Omega_k = 0$ for~$k >p$. Thus,~$i(\sum_{k=2}^p (\bar \nabla \cdot \omega_k) \omega_k)=\sum_{k=2}^n (\nabla \cdot \Omega_k) \Omega_k)$. Hence, from~\eqref{eq:om1_reduced} we obtain~\eqref{eq:om1}.

Similarly,~$i((\omega_j \cdot \bar \nabla) \omega_1)=(\Omega_j \cdot \nabla) \Omega_1.$ and~$i((\bar \nabla \omega_j) \omega_1)=(\nabla \Omega_j) \Omega_1$ (details are left to the reader), and from~\eqref{eq:omj_reduced} we obtain~\eqref{eq:omj} (for~$1\leq j \leq p$). The remaining equations~\eqref{eq:omj} are satisfied because all terms are constant (in space and time). Therefore~$\rho,\rot$ satisfies the SOHB system~\eqref{eq:mass_2}-\eqref{eq:expressF} in dimension~$n$, since it is equivalent to~\eqref{eq:mass_rhom}-\eqref{eq:omj}, and the initial conditions are exactly~$\rho^0$ and~$\rot^0$. By uniqueness, this prove that the solution is actually of the desired form, and ends the proof of Proposition~\ref{prop_dim_reduc}. We notice that this also proves existence and uniqueness of the solution to the SOHB system in dimension~$p$, on the same interval~$[0,T]$ as the solution in dimension~$n$, even if we only required local existence in time in dimension~$p$ and existence and uniqueness in dimension~$n$.
\endproof

\end{document}